\documentclass[11pt,a4paper]{amsart}
\usepackage[utf8]{inputenc}
\usepackage{amsmath}
\usepackage{amscd}
\usepackage{amsfonts}
\usepackage{amssymb}
\usepackage{amsthm}
\usepackage{caption}
\usepackage{color}
\usepackage{graphicx}
\usepackage{scrextend}                                           
\usepackage{xcolor}
\usepackage{wrapfig}          

\addtokomafont{labelinglabel}{\sffamily}

\def\R{\mathbb{R}}
\def\N{\mathbb{N}}
\def\H{\mathbb{H}^2}
\def\bL{\mathbb{L}}
\def\BL{\mathbb{L}_{\bsl}}

\def\mrd{\mathrm{d}}

\def\WP{\mathrm{WP}}

\def\n{{\mathfrak{n}}}
\def\g{\mathfrak{g}}

\def\mC{\mathcal{C}}
\def\mG{\mathcal{G}}
\def\mH{\mathcal{H}}
\def\mN{\mathcal{N}}
\def\mR{\mathcal{R}}
\def\mS{\mathcal{S}}

\def\mcML{\mathcal{ML}}
\def\mcT{\mathcal{T}}
\def\mCi{\mathcal{C}^\infty}

\def\tde{\tilde{\de}}
\def\tE{\tilde{E}}
\def\tl{{\tilde{\lambda}}}
\def\ts{\tilde{s}}
\def\tg{\tilde{\gamma}}
\def\tp{\tilde{p}}
\def\td{\tilde{\delta}}
\def\tx{\tilde{x}}

\def\ol{\overline}
\def\ul{\underline}
\def\ulx{\underline{x}}
\def\uly{\underline{y}}
\def\ulz{\underline{z}}

\def\de{\partial}
\def\dinf{\partial_\infty}

\def\psl{PSL(2,\mathbb{R})}

\def\SO{SO_0(2,1)}
\def\Mink{\mathbf{M}^3}
\def\mRb{\mR(\mbb)}
\def\Hyp{\H}
\def\piS{\pi_1(S)}

\def\mT{\mathcal{T}_S}
\def\omT{\ol{\mathcal{T}}_S}
\def\mTb{\mathcal{T}_S(\mbb)}
\def\TmTr{T\mT\restr{\mTb}}

\def\mML{\mathcal{ML}_S}
\def\mMLu{\mathcal{ML}_S^\#}
\def\mCML{\mathcal{ML}^c_S}
\def\mFMLu{\mathcal{F}\mMLu}
\def\mMLd{\mML^\dagger}

\def\pl{{\lambda_+}}
\def\ml{{\lambda_-}} 

\def\lvn{\Lambda_n^{[v]}}

\def\lcv{\Lambda^{[cv]}}

\def\lpm{{(\lambda_+,\lambda_-)}}
\def\bsl{{\boldsymbol\lambda}}
\def\bsL{{\boldsymbol\Lambda}}
\def\bsm{{\boldsymbol\mu}}
\def\hmu{{\hat{\mu}}}

\def\mbb{\mathbf{b}}
\def\S{S}

\def\ddt{\frac{\mathrm{d}}{\mathrm{d}t}}
\def\Egr{E^{t\gamma}_l}
\def\zts{(\zeta_t)_*}
\def\Psib{\Psi^{\mbb}}
\def\Psibo{\Psi_0^{\mbb}}
\def\vv{(\O,\O)}
\def\hta{{(h,\tau)}}

\DeclareMathOperator{\Hess}{Hess}

\DeclareMathOperator{\Diff}{Diff}
\DeclareMathOperator{\Cos}{Cos}
\DeclareMathOperator{\supp}{supp}
\DeclareMathOperator{\Isom}{Isom}
\DeclareMathOperator{\TL}{T}

\DeclareMathOperator{\meas}{meas}

\DeclareMathOperator{\diag}{diag}
\DeclareMathOperator{\cmp}{cmp}
\DeclareMathOperator{\Stab}{Stab}

\DeclareMathOperator{\spir}{spir}

\newcommand*\quot[2]{{^{\textstyle #1}\big/_{\textstyle #2}}}
\newcommand{\restr}[1]{\hspace{-0.3em}\mid_{#1}}
\newcommand{\e}[1]{e^{#1}_l}

\renewcommand{\emptyset}{\O}

\renewcommand{\setminus}{\smallsetminus}

\def\tsd{\theoremstyle{definition}}
\def\tsp{\theoremstyle{plain}}
\def\tsr{\theoremstyle{remark}}

\tsp
\newtheorem{thmPaper}{Theorem}

\tsp

\tsd

\tsr

\tsr

\begin{document}

\title{Hamiltonian properties of earthquake flows on surfaces with closed geodesic boundary}
\author{Daniele Rosmondi}
\address{Dipartimento di Matematica \lq Felice Casorati\rq, Universit\`a degli Studi di Pavia, Via Ferrata 5, 27100 Pavia, Italy}
\email{daniele.rosmondi01@ateneopv.it}
\thanks{Partially supported by FIRB project \lq Geometry and topology of low-dimensional manifolds\rq }

\begin{abstract}
The Teichm\"uller space $\mTb$ of hyperbolic metrics on a surface $S$ with fixed lengths at the boundary components is symplectic. We prove that any sum of infinitesimal earthquakes on $S$ that is tangent to $\mTb$ is Hamiltonian, by providing a Hamiltonian $\bL$. Such function extends the classical length map associated to a compactly supported measured geodesic lamination and shares with it some peculiar properties, such as properness and strict convexity along earthquakes paths under usual topological conditions. As an application, we prove that any non-Fuchsian affine representation of $\piS$ into $\R^{2,1}\rtimes\SO$ with cocompact discrete linear part is determined by the singularities of the two invariant regular domains in $\R^{2,1}$ pointed out by Barbot, once the boundary lengths are fixed.
\end{abstract}

\maketitle


\pagestyle{plain}

\section*{Introduction}
\label{Sez0}
Let $S$ be a surface of genus $\g$ with $\n$ closed mutually disjoint disks removed, with $\chi(S)=2-2\g-\n<0$. Consider the space $\mT$ of hyperbolic metrics on $S$ whose completion $\ol{S}$ has $\n$ closed geodesic boundary components $\de_1,\ldots,\de_\n$, up to diffeomorphisms of $\ol{S}$ isotopic to the identity. Such metrics can be deformed via \textit{left/right hyperbolic earthquakes}, which roughly speaking transform $h\in\mT$ to $h'\in\mT$ by shearing $(S,h)$ towards the left/right along \textit{measured geodesic laminations}, whose space is denoted by $\mML$.  Weighted closed geodesics are the basic examples of elements of $\mML$. Thus, associated with each measured geodesic lamination $\lambda$ there are the left and right earthquake maps $E^\lambda_l,E^\lambda_r\colon\mT\to\mT$.\\
Let us first consider when $S$ is closed, i.e. $\n=0$. The space of weighted closed geodesics is in this case dense in $\mML$. With every $\lambda\in\mML$ it is associated the \textit{length map} $L_\lambda\colon \mT\to\R$, defined for any $\omega$-weighted closed geodesic $c$ as $L_\lambda(h)=\omega\,\ell_h(c)$ and extended for $\lambda\in\mML$ by approximation. It was proved by Wolpert in \cite{W2} that {$E^\lambda_l$} is the Hamiltonian flow of $-L_\lambda$ with respect to the Weil-Petersson form $\varpi_{\WP}$ on $T\mT$. The related Hamiltonian vector field is denoted by {$e^\lambda_l$}.\\
The aim of this paper is to extend such result when $\n>0$. In such attempt, some tools and certain statements occurring in the closed case go missing. First of all, $\lambda\in\mML$ can contain geodesics \textit{spiralling} near boundary components of $S$. This implies that $\lambda$ can not be approximated by weighted closed geodesics, and a priori it is not clear how a length map $L_\lambda$ can be defined. Moreover, $\mT$ is no longer a symplectic manifold (its dimension could even be odd). This can be bypassed by partioning $\mT$ with submanifolds which are symplectic: for every $\mbb=(b_1,\ldots,b_\n)\in(\R_{>0})^\n$, on the tangent of $\mTb$, the space of metrics $h\in\mT$ with fixed boundary lengths $b_i=\ell_h(\de_i)$, a symplectic structure is induced by the one on $T\mcT_{2S}$, where $2S$ denotes the \textit{double} of $S$. However, if $\lambda\in\mML$ has spiralling leaves then the \textit{infinitesimal {left} earthquake} {$e^\lambda_l$}$\in\Gamma(T\mT)$ is not tangent to $\mTb$. {There is a notion of signed intersection of a lamination $\lambda$ near a boundary component $\de_i$ (see \cite{B-K-S}). For any $N$-uple $\bsl=(\lambda_1,\ldots,\lambda_N)\in\mML^{\ N}$, the vector field $\e{\bsl}=\e{\lambda_N}+\ldots+\e{\lambda_1}$ is tangent to $\mTb$ if and only if the sum of the signed intersections of $\lambda_1,\ldots,\lambda_N$ near $\de_i$ is null for every $i$. We denote the space of such $N$-uples by $\mMLu$.} The main theorem can now be stated as follows.

\begin{thmPaper}
\label{Mainthm}
Given $\mbb\in(\R_{>0})^\n$, the vector field {$\e{\bsl}$}$\in\Gamma(T\mTb)$ is Hamiltonian for every $\bsl\in\mMLu$.
\end{thmPaper}

\noindent
We provide a Hamiltonian $-\BL\colon \mTb\to\R$ which extends $-\sum L_{\lambda_n}$ to the case when $\n>0$. We also show that $\BL$ is strictly convex (in a suitable sense) and proper if $\bsl\in\mMLu$ is a $N$-uple that \textit{fills up} $S$, i.e. every simple closed non-trivial and non-peripheral curve meets the support of $\lambda_1\cup\ldots\cup\lambda_N$. We denote by $\mFMLu$ the space of filling couples in $\mMLu$.\\
{We provide an application within the study of flat Lorentzian structures, analogue to the compact case shown in \cite{B-S}. Identifying $\R^{2,1}$ with the Lie algebra of $SL(2,\R)$ (through the Killing form) and $\mT$ with the space $\mR$ of Fuchsian cocompact representations of $\piS$, the tangent space $T\mT$ can be identified with the space of affine deformations of elements of $\mR$. Barbot showed in \cite{Ba3} that associated with $\rho=h+\tau\in\mR$ there are two $\rho$-regular domains (as they are called in \cite{B-B}) in $\R^{2,1}$. Each domain is determined by a lamination on the surface base point of $\rho$, viewed as the dual of the singularities of the domains (see \cite{B-B}). The couple $\lpm$ of such laminations fills up $S$ and satisfies the condition $\tau=\e{\ml}(h)=-\e{\pl}(h)$. We show that $\lpm$ determines $\rho$ up to fixing the boundary lengths:

\begin{thmPaper}
\label{PaperThmB}
The map $\Psi:T\mT\to\mML^{\ 2}$ associating $\rho=(h,\tau)$ with the couple $\lpm$ described above is a fibration over $\mFMLu$, the subset of $\mMLu$ of filling couples. The fiber is isomorphic to $\R^\n$.
\end{thmPaper}
}
This paper is organized as follows. In the first part of Section \ref{Sez1} we recall general notions about measured geodesic laminations and hyperbolic earthquakes on $S$. After that, we proceed to give to $\mML$ a manifold structure compatible with the weak$^*$-convergence topology and we study smoothness of infinitesimal earthquakes. Finally, we endow $\mTb$ with a symplectic structure.\\
Section \ref{Sez2} is devoted to the construction of $\BL$, starting from the Hamiltonian condition and decomposing any $\bsl\in\mMLu$ in the union of simple couples in $\mMLu$, in a suitable sense. After defining $\bL$ for these simple couples and checking the Hamiltonian condition, we provide $\BL$ for generic $\bsl\in\mMLu$. Properness and strict convexity of $\BL$ are proved in Section \ref{Sez3}, where is also computed $\Hess\BL$ at its critical point.\\
In Section \ref{Sez4} we apply such results to the study of $\Psi\colon T\mT\to\mFMLu$.


\section{Earthquakes and measured geodesic laminations}
\label{Sez1}

Given a topological surface $S$ obtained by removing $\n$ closed mutually disjoint disks from a compact surface of genus $\g$ with Euler characteristic $\chi(S)=2-2\g-\n<0$, let
\begin{align*}
\mT=\{&\text{hyperbolic metrics on $S$ whose completion $\ol{S}$}\\
&\text{has $\n$ closed geodesic boundary components}\}/\Diff_0(\ol{S}),
\end{align*}
where $\Diff_0(\ol{S})$ denotes the group of the diffeomorphisms of $\ol{S}$ isotopic to the identity. We will refer to the boundary components of $S$ as $\de_1,\ldots,\de_\n$.

\subsection{Measured geodesic laminations}
\label{S_ML}

\tsd
\newtheorem{meas.lam}[d1]{Definition}
\begin{meas.lam}
Given a hyperbolic metric $h$ on $S$, a geodesic lamination on $(S,h)$ is the data $\lambda$ of a family of mutually disjoint complete simple geodesics (called the leaves of $\lambda$) whose union is a closed subset (called the support of $\lambda$ and denoted by $\supp(\lambda)$) of $S$. 
A measured geodesic lamination of $S$ is the data of a geodesic lamination $\lambda$ and a transverse measure $\meas_\lambda$, that is a measure defined on the arcs on $S$ transverse to each leaf of $\lambda$ and  with endpoints in $S\setminus\supp(\lambda)$ such that 
\begin{itemize}
\item[-] $\meas_\lambda(c)\neq 0$ if and only if $c\cap\supp(\lambda)\neq \emptyset$;
\item[-] if there exists an isotopy between two arcs $c_1$ and $c_2$ realized through arcs transverse to $\lambda$ then $\meas_\lambda(c_1)=\meas_\lambda(c_2)$.
\end{itemize}  
\end{meas.lam}

Weighted multicurves are the simplest examples of measured geodesic lamination on $\S$. The support is the finite union of simple closed mutually disjoint non trivial geodesics $\gamma_i$. Chosen real positive numbers $\omega_i$ (called weights) respectively assigned to $\gamma_i$, the transverse measure is given by
\[
c\mapsto\sum_i\omega_i\cdot\#(\gamma_i\cap c)
\]
for any arc $c$ transverse to $\bigcup\gamma_i$.\\
It is known (see \cite{C}) that the Lebesgue measure of the support of a geodesic lamination is zero.\\
\begin{wrapfigure}{r}{0.35\textwidth}
\centering
\includegraphics[width=0.9\linewidth]{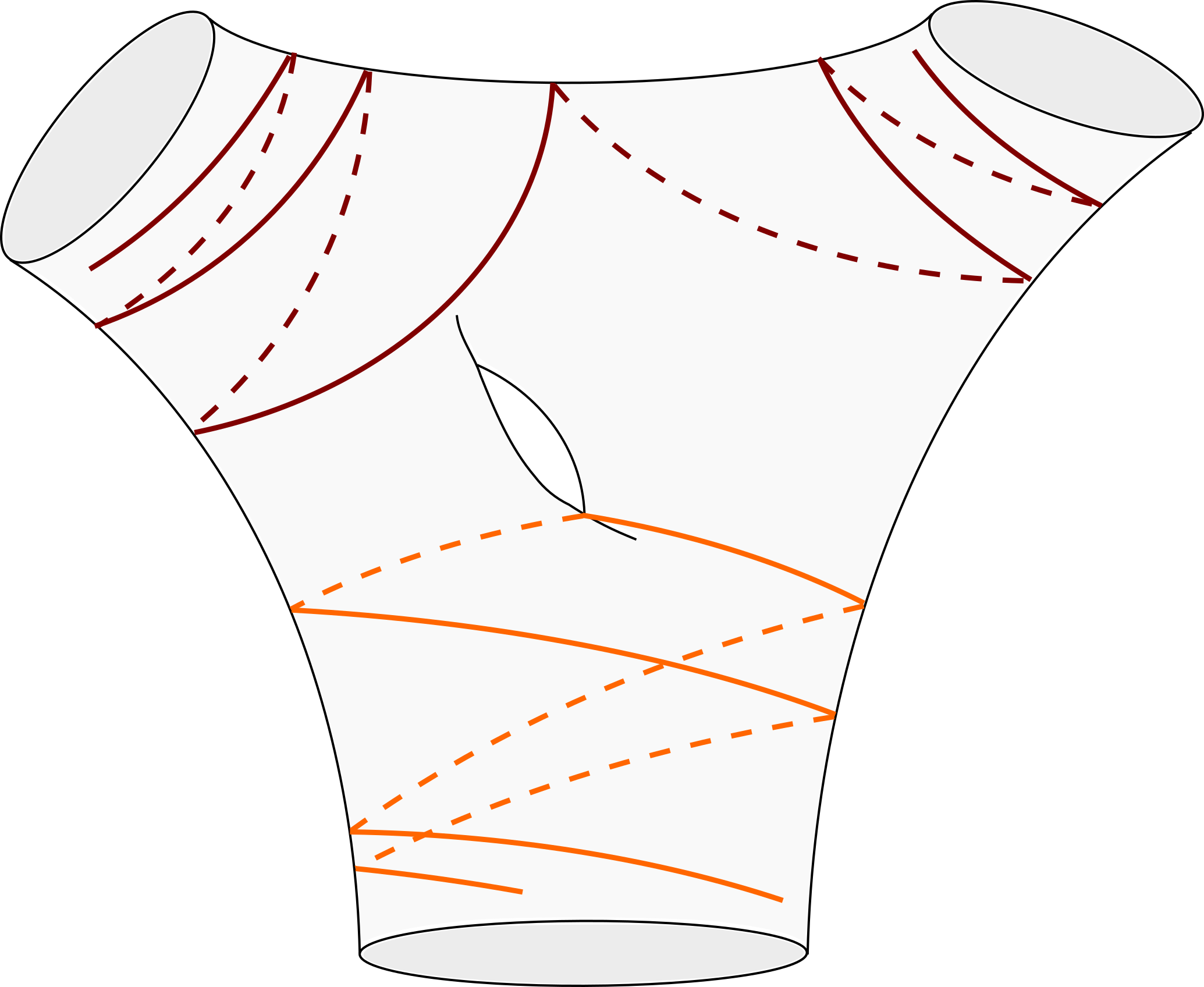} 
\caption{A geodesic lamination with two spiralling leaves}
\label{Spyr_Lam}
\end{wrapfigure}If $h\in\mT$ then any measured geodesic lamination $\lambda$ on $(S,h)$ has a maximal compact sublamination $\lambda^{(0)}$, in the sense that if $\mu$ is a sublamination of $\lambda$ with compact support in $S$ then $\mu$ is a sublamination of $\lambda^{(0)}$ too. Each leaf of $\supp(\lambda)\setminus\supp(\lambda^{(0)})$ is homeomorphic to $\R$ and spirals near two boundary components (possibly coincident) of $S$ (see Figure \ref{Spyr_Lam}).\\
%
If we denote by $\mcML_{(\S,h)}$ the measured geodesic laminations on $(\S,h)$ with $h\in\mcT_\S$, being a space of measures it seems natural to provide it with the topology of the weak-convergence of measures (sometimes also called weak$^*$-convergence). It is known (see Section 1.7 of  \cite{P-H}) that for every $h_1,h_2$ in $\mcT_\S$ there is a homeomorphism $F:\mcML_{(\S,h_1)}\to\mcML_{(\S,h_2)}$ so that, roughly speaking, $\supp(F(\lambda))$ is obtained straightening with respect to $h_2$ the leaves of $\supp(\lambda)$. This suggests that it makes sense to associate $\mcT_\S$ with the space $\mcML_\S$ of measured laminations, whose support is only a topological data; this space inherits the weak convergence topology. Finally, define
\[
\mCML=\{\lambda\in\mcML_S\,|\,\lambda=\lambda^{(0)}\},
\]
the subspace of laminations with compact support. The following theorem is a well known result (see \cite{P-H}).
\begin{figure}[h]
\centering
\includegraphics[scale=.15]{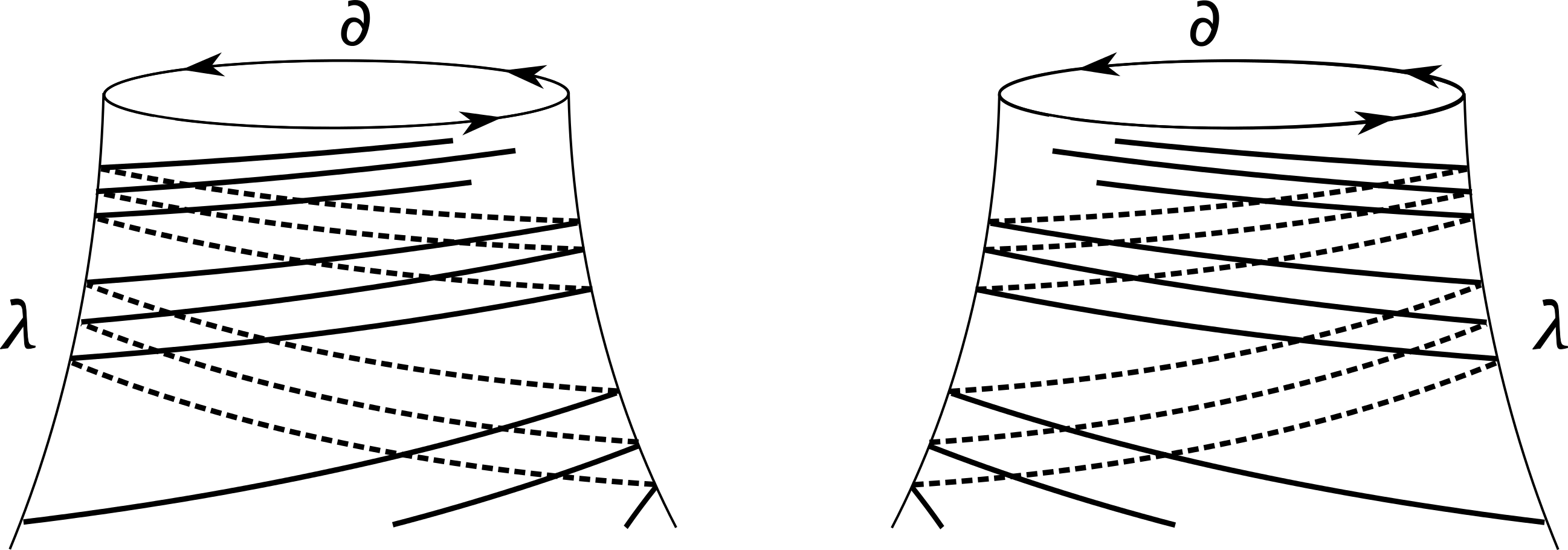}
\caption{Respectively, positive and negative sense of spiralling}
\label{Spyr_sense}
\end{figure}
\tsp
\newtheorem{multic.dense}[p1]{Theorem}
\begin{multic.dense}
The space of weighted multicurves on $\S$ is dense in $\mCML$.
\end{multic.dense}

\begin{wrapfigure}{r}{0.45\textwidth}
\centering
\includegraphics[width=1.0\linewidth]{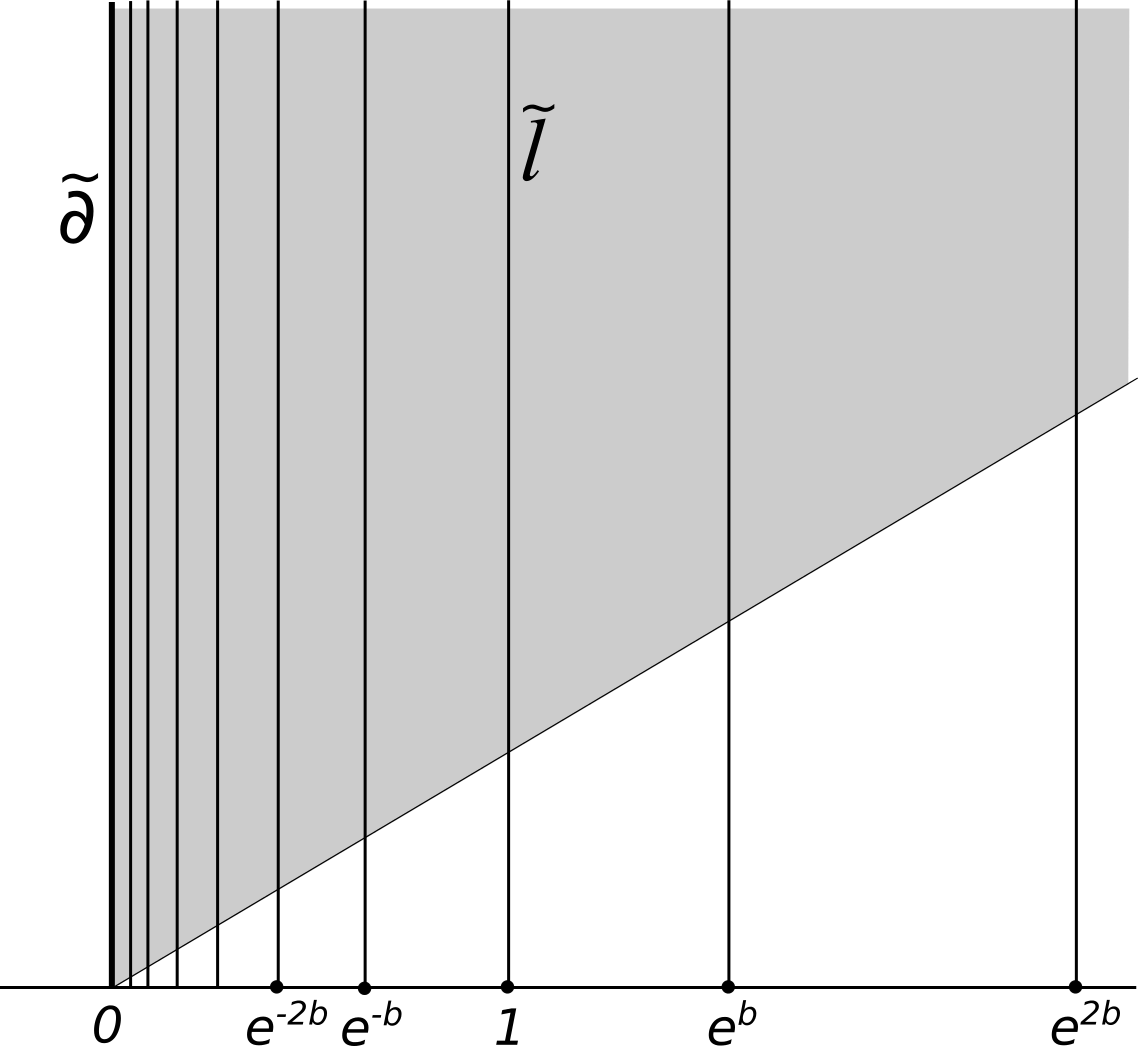}
\caption{}
\label{N_eps_fig}
\end{wrapfigure}
Let us fix for a moment $h\in\mT$ and consider a measured geodesic lamination $\lambda$ on $(S,h)$. If a leaf of $\lambda$ is not contained contained in a compact subset of $S$, then, in order to be a complete geodesic with no self-intersections, it must spiral along one or two connected components of $\de S$. There are two possible senses of spiralization, as shown in Figure \ref{Spyr_sense}.\\
In particular, 
if a leaf $l$ of $\lambda$ spirals near $\de$, then for every lift $\tde\subset\H$ of $\de$ there is an $\varepsilon$-neighbourhood of $\tde$ where the preimage of $l$ is the $\Stab(\tde)$-orbit of any lift $\tilde{l}$ of $l$ sharing an ideal endpoint of $\tde$, as in Figure \ref{N_eps_fig}. See also Lemma \ref{B.collar}.\\
\\
It is possible to define the \textit{mass} $\iota(\de,\lambda)$ of $\de$ with respect to $\lambda$, a positive number that encodes how much the measure of $\lambda$ is concentrated near $\de$. It is constructed as follows. For every $x\in N_\varepsilon(\de)$ denote by $c_x$ the loop with vertex at $x$ parallel at $\de$ such that $c_x\setminus\{x\}$ is an open geodesic arc. Since $\meas_\lambda(c_x)=\meas_\lambda(c_y)$ for every $x,y\in N_\varepsilon$, as shown in \cite{B-S}, Subsection 2.3, it is well defined the mass $\iota(\de,\lambda)=\meas_\lambda(c_x)$. Moreover, $\iota(\de,\lambda)=0$ if and only if $\supp(\lambda)\cap N_\varepsilon=\O$. The mass of $\de$ does not take in account in which sense $\lambda$ spirals. Fix once for all an orientation of $\de\S$. Such choice defines a positive and a negative sense of spiralization around $\de$, as in Figure \ref{Spyr_sense}. It is now possible to define the \textit{signed mass} $m(\de,\lambda)$ of $\de$ with respect to $\lambda$ as
\begin{equation}
\label{Sign_mass}
m(\de,\lambda)=\left\{\begin{array}{rl}
+\iota(\de,\lambda) \text{if $\lambda$ spirals in the positive sense around $\de$}\\
-\iota(\de,\lambda) \text{if $\lambda$ spirals in the negative sense around $\de$}
\end{array}\right. .
\end{equation}

\tsr
\newtheorem{tigo}[r1]{Remark}
\begin{tigo}
The signed mass of $\de$ with respect to $\lambda$ is positive (respectively negative) if and only if for every orientated lift of $\de$ on $\mH$ its ending (respectively starting) ideal endpoint is contained in the set of the ideal points of the whole preimage of $\lambda$.
\end{tigo}

\subsection{Hyperbolic earthquakes}
\label{S_quakes}

Let $\mH$ be a convex subset of $\H$ with geodesic boundary.

\tsd
\newtheorem{def.earthq.surf}[d1]{Definition}
\begin{def.earthq.surf}
Given a geodesic lamination $\lambda$ in $\mH$, a left (respectively right) hyperbolic earthquake on $\mH$ along $\lambda$ is an injective (possibly discontinuous) map $\tE:\mH\to\H$ such that \\
\begin{itemize}
\item[-] the restriction of $\tE$ on a stratum of $\lambda$ is an isometry;
\item[-] denoting by $A_F\in PSL(2,\R)$ the isometry of $\H$ extending $\tE\restr{F}$ for every stratum $F$, the \textit{comparison map}
\[
\cmp(F,G)=A_F^{-1}\circ A_G:\H\to\H
\]
between two different strata $F$ and $G$ of $\lambda$ is a hyperbolic transformation whose axis weakly separates $F$ and $G$ and which translates to the left (respectively right), as viewed from $F$.
\end{itemize}
The lamination $\lambda$ is called \textit{fault locus} of the earthquake $\tE$.\\
It turns out that $\tE(\mH)$ is still a convex subset of $\H$ with geodesic boundary, as a consequence of Lemma 8.4 in \cite{B-K-S}.\\
Given a surface $\S$ and two hyperbolic metrics $h_1,h_2$ on $\S$, set $\S_i=(\S,h_i)$ for $i=1,2$. Suppose that the universal covering $\mH_i\subset\H$ of $\S_i$ is convex with geodesic boundary. A bijective map $E:\S_1\to\S_2$ is a left (respectively right) hyperbolic earthquake if it has a lifting $\tE:\mH_1\to\mH_2$ which is a left (respectively right) hyperbolic earthquake on $\mH_1$.
\end{def.earthq.surf}

%
The fault locus can be endowed with a transverse measure encoding the shearing of the earthquake, obtaining a measured geodesic lamination: the $\omega$-weighted curve $c$. This can be done in general, as stated in the following (\cite{T}, Proposition 6.1).

\tsp
\newtheorem{met.earthq}[p1]{Proposition}
\begin{met.earthq}
A measured geodesic lamination $\lambda\subset\mH$ is associated to any earthquake so that $\supp(\lambda)$ coincides with the fault locus; if $a:[0,1]\to\mH$ is an arc with endpoints in $\mH\setminus\lambda$ then
\[
\meas_\lambda(a)=\inf_{P\text{ partition of [0,1]}}\sum_{i=1}^{I_P}\TL(\cmp(A_{F_{i-1}},A_{F_i}))
\]
where for every partition $P=(0=t_0,t_1,t_2,\ldots,t_{I_P}=1)$ of $[0,1]$ the stratum $F_i$ of $\lambda$ is the one containing $t_i$. Here $\TL(B)$ denotes the translation length of a hyperbolic transformation $B$.\qed
\end{met.earthq}

Moreover, Thurston showed that different earthquakes produce different measured geodesic laminations (see \cite{T}). The converse holds, since we did not suppose that $\tE$ is surjective. See \cite{B-K-S} for further details.\\
%
There is a natural non surjective immersion of $\mT$ into the Teichm\"uller space $\omT$ of hyperbolic metrics on $S$ of finite area whose completion has compact geodesic boundary. A metric in $\omT$ can have cusps at some punctures of $S$. %
%
Associated with $\lambda\in\mML$, there are a left and a right earthquake map between $\mT$ and $\omT$: 
\[
E^\lambda_l,E^\lambda_r\colon \mT\to\omT.
\]
Proposition 3.3 in \cite{B-K-S} shows explicitly how right and left earthquakes change the length of the boundary components $\de_1,\ldots,\de_\n$ of $S$: for every $h\in\mT$ and $\lambda,\mu\in\mML$
\begin{equation}
\label{M}
\left\{\begin{array}{rl}
\ell_{E^\lambda_l(h)}(\de_i)=|\ell_h(\de_i)-m(\de_i,\lambda)|\\
\ell_{E^\mu_r(h)}(\de_i)=|\ell_h(\de_i)+m(\de_i,\mu)|
\end{array}\right. .
\end{equation}
Fix $\mbb=(b_1,\ldots,b_\n)\in(\R_{>0})^\n$ and set
\begin{equation*}
\mTb=\{h\in\mT\,|\,\ell_h(\de_i)=b_i\ \ \forall i=1,\ldots,\n\}.
\end{equation*}
Clearly, 
\[
\mT=\bigcup_{\mbb\in(\R_{>0})^\n}\mTb.
\]
{In this paper we are interested in $N$-uples $\bsl\in\mML^{\ N}$ for which the vector field 
\[
\e{\bsl}(h)=\ddt_{|0}\Big(E^{t\lambda_N}_l\circ\ldots\circ E^{t\lambda_1}_l(h)\Big)
\]
is tangent to $\mTb$, with $\mbb\in(\R_{>0})^\n$. Now, for every $i=1,\ldots,\n$, if $h_t=E^{t\lambda_N}_l\circ\ldots\circ E^{t\lambda_1}_l(h)$ then, using \eqref{M} for $t$ sufficiently small,
\[
\ell_{h_t}(\de_i)=\ell_{h}(\de_i)-tm(\de_i,\lambda_1)-tm(\de_i,\lambda_2)-\ldots-tm(\de_i,\lambda_N)
\]
and so $\e{\bsl}\in\Gamma(T\mTb)$ if and only if 
\[
0=\ddt_{|0}\ell_{h_t}(\de_i)=-m(\de_i,\lambda_1)-m(\de_i,\lambda_2)-\ldots-m(\de_i,\lambda_N)
\]
for every $h\in\mTb$ and $i=1,\ldots,\n$. Notice that such condition is actually independent on $\mbb$. Thus, we introduce the space
\[
\mMLu=\bigg\{\bsl\in\mML^{\ N}\,\bigg|\,\sum_{n=1}^Nm(\de_i,\lambda_n)=0\ \  \forall i=1,\ldots,\n\bigg\}.
\]
\tsr
\newtheorem{N1case}[r1]{Remark}
\begin{N1case}
When $N=1$ then $\mMLu=\mCML$. Since classical results are already known for compactly supported laminations, we will suppose from now on that $N\geq 2$.
\end{N1case}
}

\subsection{The topology of $\mML$}
\label{S_topML}

{Now we are going to give to $\mML$ a manifold structure. First let us introduce the straightening $\nu^R$ of a measured lamination $\nu\in\mML$. If $\gamma$ is a spiralling geodesic between two connected components $\de_i$ and $\de_j$ of $\de \S$, consider its preimage $\Gamma$ on the universal cover $\mH\subset \H$. Every connected component of $\Gamma$ is a geodesic $\tilde{\gamma}$ with endpoints in the (ideal closure) of certain lifts $\tde_i$ and $\tde_j$ of $\de_i$ and $\de_j$ respectively. If we replace each $\tilde{\gamma}$ with the geodesic arc $\tilde{\gamma}^R$ with endpoints on $\tde_i$ and $\tde_j$ perpendicular to $\tde_i$ and $\tde_j$ and we project $\tilde{\gamma}^R$ on $\S$, we obtain a geodesic arc $\gamma^R$ on $\S$ normal to $\de_i$ and $\de_j$ with endpoints on $\de_i$ and $\de_j$. For each $\nu\in\mML$ denote by $\nu^R$ the set of geodesic (weighted) arcs obtained by $\nu$ replacing each spiralling geodesic $\gamma$ of $\nu$ with $\gamma^R$.\\
Consider the set $\mML^R=\{\nu^R\,|\,\nu\in\mML\}$. This space is a submanifold of the space of measured laminations (that we denote by $\mMLd$) studied in \cite{A-L-P-S}; we will mention only the necessary details. Using the notation of \cite{A-L-P-S}, we fix a pant decomposition 
\[
P=\{C_1,\ldots,C_{3\g-3+\n},B_1=\de_1,\ldots,B_\n=\de_\n\}
\]
of $\S$ with internal curves $C_1,\ldots,C_{3\g-3+\n}$ and boundary curves $B_1=\de_1,\ldots,$ $B_\n=\de_\n$. Every lamination $\sigma\in\mMLd$ has coordinates  
\[
\big(DT(\sigma,C_1),\ldots,DT(\sigma,C_{3\g-3+\n}),\hat{\theta}(\sigma,B_1),\ldots,\hat{\theta}(\sigma,B_\n)\big)
\]
where $DT(\sigma,C_i)\in\R^2$ depends on the behaviour of $\sigma$ in a regular neighbourhood of $C_i$ and $\hat{\theta}(\sigma,\de_i)\in\R$ depends on the behaviour with respect to the boundary component $\de_i$. Following their constructions, it turns out that, for every $\nu\in\mML$, $\hat{\theta}(\nu^R, \de_i)=\iota(\nu,\de_i)\geq 0$. So if we consider the coordinates $\Theta_P\colon \mathcal{ML}\to\R^{6\g-6+3\n}$ such that
\begin{equation}
\label{Thetacoord}
\Theta_P(\nu)=\big(DT(\nu^R,C_1),\ldots,DT(\nu^R,C_{3\g-3+\n}),m(\nu,\de_1),\ldots,m(\nu,\de_\n)\big)
\end{equation}
for $\nu\in\mML$, where $m(\nu,\de_i)$ is the signed mass defined by \eqref{Sign_mass}, we provide $\mML$ with a manifold structure. Such coordinates depend on the pant decomposition $P$; however, if $P'$ is another pant decomposition, notice that the last $\n$ coordinates does not depend on the pant decomposition, whereas applying the results in \cite{A-L-P-S} the change of coordinates of the other components is smooth.\\
Even if the projection $\mML\to\mML^R$ is not injective, the map $\Theta_P$ is injective, since we have avoided the ambiguity given by the spiralling senses around $\de\S$. \\
It is shown in \cite{A-L-P-S} that the topology on $\mMLd$ coincides with the topology of the weak$^*$-convergence of measures. We are interested to show that also for $\mML$ the topology is the one of weak$^*$-convergence of measures.
\tsp
\newtheorem{wsconv_lem}[p1]{Lemma}
\begin{wsconv_lem}
\label{Wsc}
Consider a sequence $\lambda_n$ converging to $\lambda$ in the manifold $\mML$. If $\lambda^{[s]}$ is the sublamination of $\lambda$ made by spiralling leaves, then the support of $\lambda^{[s]}$ is contained in $\lambda_n$ for $n$ sufficiently big. In particular, there exist decompositions 
\begin{gather*}
\lambda_n=\lambda_n^{[c]}\cup\lambda_n^{[s]}\cup\lambda_n^{[v]},\\
\lambda=\lambda^{[cc]}\cup\lambda^{[s]}\cup\lambda^{[cv]}
\end{gather*}
such that, up to passing to a subsequence,
\begin{itemize}
\item $\lambda_n^{[c]}$ is the maximal compact sublamination of $\lambda_n$, and $\lambda_n^{[c]}$ converges to $\lambda^{[cc]}$;
\item $\lambda^{[s]}$ is the sublamination of $\lambda$ whose support consists of the spiralling leaves of $\lambda$, and $\lambda_n^{[s]}$ is the maximal sublamination of $\lambda_n$ such that $\supp(\lambda_n^{[s]})=\supp(\lambda^{[s]})$; moreover, $\lambda_n^{[s]}$ tends to $\lambda^{[s]}$;
\item $\lambda_n^{[v]}$ is the complementary of $\lambda_n^{[s]}$ in the spiralling part of $\lambda_n$, so that $\lambda_n^{[v]}$ converges to the compact lamination $\lambda^{[cv]}$.
\end{itemize}
\begin{proof}
{We prove that if $l_n$ is a sequence of leaves of $\lambda_n$ converging to a leaf $l\in\lambda^{[s]}$, then $l_n=l$ for $n$ big. The claim directly implies the statement. Let us prove the claim.\\
Consider a leaf $l$ of $\lambda^{[s]}$, going say between the boundary components $\de$ and $\de'$ of $\S$. On the universal covering $\mH\subset\H$ of $\S$, consider a lift $\tilde{l}$ of $l$, going from $\tde$ and $\tde'$, the boundary components of $\de\mH$ who projects onto $\de$ and $\de'$ respectively. The straightening $\tilde{l}^R$ of $\tilde{l}$ has an endpoint $z\in\tde$. There is a $\delta$-neighbourhood $U$ of $\tilde{l}^R$ in $\mH$ such that for every $u\in(\ol{U}\cap\tde)\setminus\{z\}$ the complete geodesic of $\H$ normal to $\tde$ passing through $u$ must intersect $\tde'$, but this intersection cannot be orthogonal, so if a lamination $\nu\in \mMLd$ meets $U\cap\tde$, then it must contain the leaf $l$. Thus, leaves of $(\lambda^{[s]})^R$ must be contained in $(\lambda_n\setminus\lambda^{[c]}_n)^R$ for big $n$, and in fact $(\lambda^{[s]})^R$ must be the limit of the sublamination $(\lambda^{[ss]}_n)^R$ made by the leaves of $(\lambda_n\setminus\lambda^{[c]}_n)^R$ whose weight is not tending to zero.} 
\end{proof}
\end{wsconv_lem}

\tsp
\newtheorem{wsconv}[p1]{Proposition}
\begin{wsconv}
If $\lambda_n\to\lambda$ in $\mML$ then for every arc $\alpha$ on $\S$ with endpoints in $\S\setminus\big(\supp(\lambda)\cup\bigcup\supp(\lambda_n)\big)$ and for every $\varphi\in C^\infty_c(\alpha)$ 
\[
\int_\alpha\varphi\,\mrd(\meas_{\lambda_n})\xrightarrow{n\to\infty}\int_\alpha\varphi\,\mrd(\meas_{\lambda}).
\]
\begin{proof}
From now on, for simplicity we will write $\mrd\lambda_n$ and $\mrd\lambda$ respectively for $\mrd(\meas_{\lambda_n})$ and $\mrd(\meas_{\lambda})$.\\
Take the decomposition
\begin{gather*}
\lambda_n=\lambda_n^{[c]}\cup\lambda_n^{[s]}\cup\lambda_n^{[v]},\\
\lambda=\lambda^{[cc]}\cup\lambda^{[s]}\cup\lambda^{[cv]}
\end{gather*}
provided by Lemma \ref{Wsc}, and consider the induced decomposition on the double straightenings $\Lambda_n$, $\Lambda$ of $\lambda_n$, $\lambda$ respectively:
\begin{gather*}
\Lambda_n=\Lambda_n^{[c]}\cup\Lambda_n^{[s]}\cup\Lambda_n^{[v]},\\
\Lambda=\Lambda^{[cc]}\cup\Lambda^{[s]}\cup\Lambda^{[cv]}.
\end{gather*}
Notice that the weights of the leaves of $\lvn$ are going to 0, since the masses of $\lvn$ at the boundary of $S$ are vanishing.\\
Fixed $\epsilon>0$ and denoting by
\begin{gather*}
T_1=\bigg|\int_\alpha\varphi\,\mrd\lambda_n^{[c]}-\int_\alpha\varphi\,\mrd\lambda^{[cc]}\bigg|\\
T_2=\bigg|\int_\alpha\varphi\,\mrd\lambda_n^{[s]}-\int_\alpha\varphi\,\mrd\lambda^{[s]}\bigg|\\
T_3=\bigg|\int_\alpha\varphi\,\mrd\lambda_n^{[v]}-\int_\alpha\varphi\,\mrd\lambda^{[cv]}\bigg|
\end{gather*}
it suffices to show that for $n$ sufficiently large $T_1+T_2+T_3\leq 6\epsilon$.\\
{It is easy to estimate $T_1\leq\epsilon$ and $T_2\leq\epsilon$ for $n$ large enough, due respectively to the compact and discrete nature of the involved sublaminations. }The term $T_3$ requires more attention. First of all, let us split is as
\begin{align*}
T_3 & \leq \bigg|\int_\alpha\varphi\,\mrd\lambda_n^{[v]}-\int_{ {\alpha}} {\varphi}\,\mrd\lvn\bigg|+\bigg|\int_{ {\alpha}} {\varphi}\,\mrd\lvn-\int_\alpha\varphi\,\mrd\lambda^{[cv]}\bigg|=\\
&=\bigg|\int_\alpha\varphi\,\mrd\lambda_n^{[v]}-\int_{ {\alpha}} {\varphi}\,\mrd\lvn\bigg|+\bigg|\int_{ {\alpha}} {\varphi}\,\mrd\lvn-\int_{ {\alpha}} {\varphi}\,\mrd\Lambda^{[cv]}\bigg|.
\end{align*}
The second term of the last member is not greater then $\epsilon$ for $n$ large enough, since $\lvn\to\lcv$. Let us consider the first one. 
Fix a lift $\tilde{\alpha}$ of $\alpha$ in the universal covering of $\S$. For every leaf $\tilde{\delta}$ of the preimage of a leaf $\delta$ of $\lvn$ denote by $D_{\tilde{\alpha}}(\tilde{\delta})$ the minimum between the lengths of the two connected components of $\tilde{\delta}^R\setminus\tilde{\alpha}$ if $\tilde{\delta}^R\cap\tilde{\alpha}$ is non empty. See also Figure \ref{Dad}. There is a constant $M=M(\alpha,\epsilon)>0$ such that if $D_{\tilde{\alpha}}(\tilde{\delta})>M$ then the \begin{figure}[h]
\centering
\includegraphics[scale=.15]{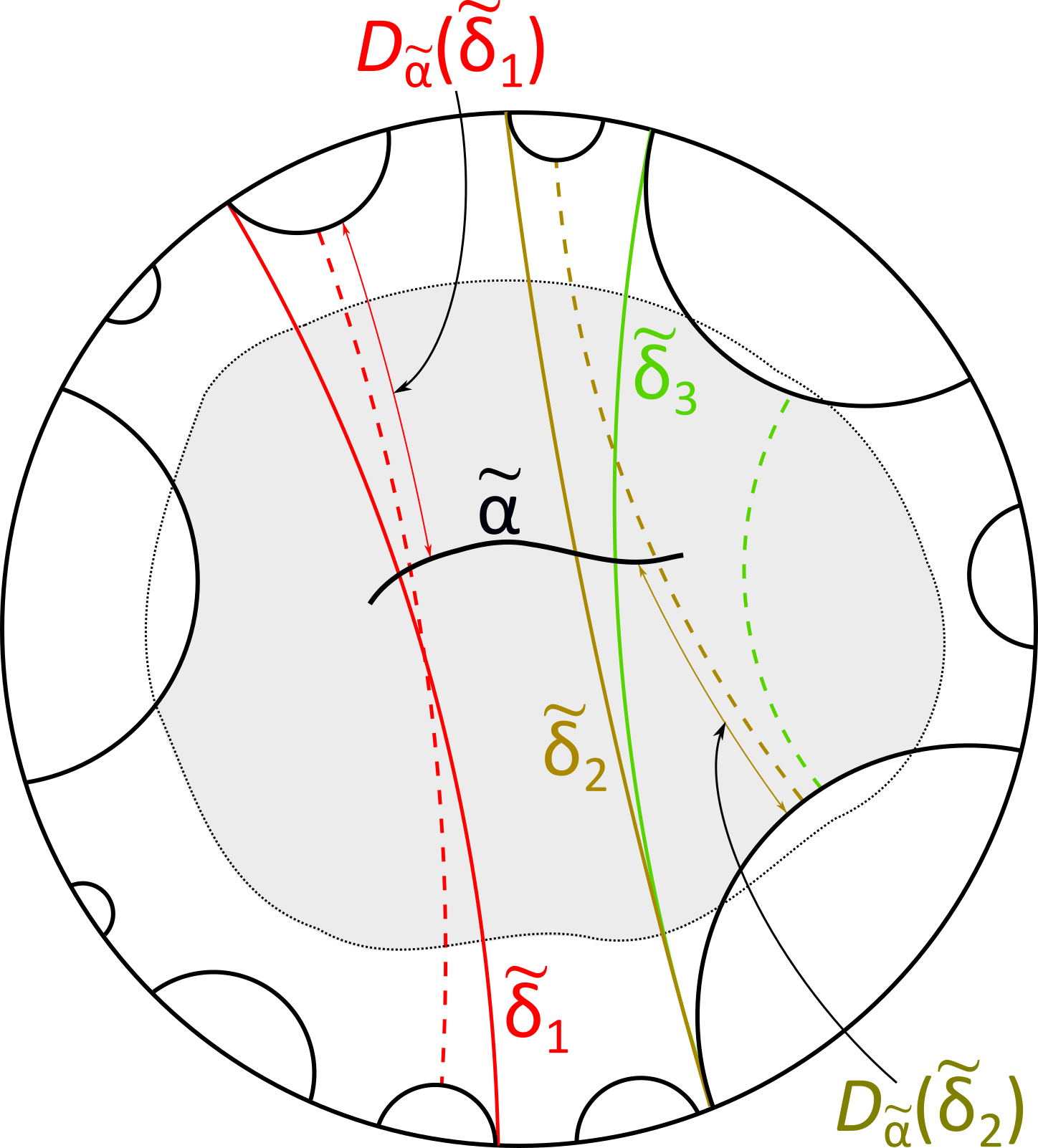}
\caption{The points in the grey region have distance from $\tilde{\alpha}$ less than $M(\alpha,\epsilon)$; the leaf $\td_1$ of $\tl^{[v]}_n$ is contained in $\tl^{[v]+}_n$, while $\td_2$ and $\td_3$ are contained in $\tl^{[v]-}_n$}
\label{Dad}
\end{figure}ideal endpoints of $\tilde{\delta}$ are close to the ones of the prolongation of $\tilde{\delta}^R$, in the Euclidean sense, so that
\[
\bigg|\int_\alpha\varphi\,\mrd\lambda_n^{[v]+}-\int_{ {\alpha}} {\varphi}\,\mrd\Lambda_n^{[v]+}\bigg|\leq \epsilon
\]
for $n$ sufficiently large, where $\lambda_n^{[v]+}$ is the sublamination of $\lambda_n^{[v]}$ of the leaves $\delta$ whose straightening meets $\alpha$ having $D_{\tilde{\alpha}}(\tilde{\delta})>M$, while $\Lambda_n^{[v]+}$ is the doubled straightening of $\lambda_n^{[v]+}$. Set $\lambda_n^{[v]-}=\lambda_n^{[v]}\setminus\lambda_n^{[v]+}$ and $\Lambda_n^{[v]-}=\Lambda_n^{[v]}\setminus\Lambda_n^{[v]+}$. Now
\begin{align*}
&\bigg|\int_\alpha\varphi\,\mrd\lambda_n^{[v]}-\int_{ {\alpha}} {\varphi}\,\mrd\lvn\bigg|\leq  \bigg|\int_\alpha\varphi\,\mrd\lambda_n^{[v]+}-\int_{ {\alpha}} {\varphi}\,\mrd\Lambda_n^{[v]+}\bigg|+\\
+&\bigg|\int_\alpha\varphi\,\mrd\lambda_n^{[v]-}-\int_{ {\alpha}} {\varphi}\,\mrd\Lambda_n^{[v]-}\bigg|\leq\epsilon + \bigg|\int_\alpha\varphi\,\mrd\lambda_n^{[v]-}\bigg|+\bigg|\int_{ {\alpha}} {\varphi}\,\mrd\Lambda_n^{[v]-}\bigg|.
\end{align*}
Actually, $\Lambda_n^{[v]-}$ (and consequently $\lambda_n^{[v]-}$) is vanishing, since its number of leaves is bounded from above by a constant depending only on the geometry of $\S$: on its universal covering $\mH$, it is easy to see that the number of connected components of $\de\mH$ distant at most $M$ from $\tilde{\alpha}$, which has compact support, are finite. Moreover, the weights of the leaves of $\lvn$ are going to 0, as $\lambda_n^{[v]}$ converges to a compact lamination. Thus, for $n$ big,
\[
\bigg|\int_\alpha\varphi\,\mrd\lambda_n^{[v]-}\bigg|+\bigg|\int_{ {\alpha}} {\varphi}\,\mrd\Lambda_n^{[v]-}\bigg|\leq 2\varepsilon.
\]
\end{proof}
\end{wsconv}
}

\subsection{Infinitesimal earthquakes}
\label{Inf.earthq}

Associated with $\lambda\in\mML$, there is the vector field
\begin{align*}
e_l^\lambda\colon &\mT\to T\mT\\
&h\mapsto\ddt_{|0}(E_l^{t\lambda}(h))
\end{align*}
called the \textit{infinitesimal left earthquake} along $\lambda$.

\tsp
\newtheorem{e_liscia}[p1]{Proposition}
\begin{e_liscia}
\label{E_liscia}
For every $\lambda\in\mML$, the vector field $e^\lambda_l$ is a smooth vector field on $\mT$.
\begin{proof}
Let us suppose $\lambda$ has a non empty compact sublamination. Decompose $\lambda=\lambda_c\oplus\lambda_s$ as the sum of the compact maximal sublamination with the spiralling sublamination. Then $e^l_l$ can be decomposed as $e_l^{\lambda_c}+e_l^{\lambda_s}$. By classical results, $e_l^{\lambda_c}$ is smooth. So we can suppose $\lambda=\lambda_s$ and consider only this case.\\
It is convenient to see $\mTb$ as the space of faithful discrete representations $ h \colon \pi_1(\S)\to\psl$ with conditions that fix the images of peripheral loops, up to conjugacy. For every $ h \in\mTb$, consider the universal covering $\mH$ of $\S$ such that $ h (\pi_1(\S))\backslash\mH\cong\S$ and fix a point $z\in\mH$; the infinitesimal earthquake regarded as an element of the cohomology $H^1(\pi_1(\S), \R^{2,1})$ is represented (see \cite{M}, \cite{nM}, \cite{B-S}) by the element $e^\lambda_l( h )\colon \pi_1(\S)\to\mathfrak{so}(2,1)\cong\R^{1,2}$ has the form
\[
\gamma\mapsto\int_{\mG}v(r)\chi_{\mG(\gamma)}(r)\,\mrd\lambda
\]
where 
\begin{itemize}
\item the space $\mathfrak{so}(2,1)$ is the Lie algebra of $SO(2,1)$,
\item the space
\[
\mG\cong (S^1\times S^1)\setminus\diag(S^1)
\]
is the set of oriented geodesics on $\H$, 
\item the map 
\[
v\colon \mG\to\mathfrak{so}(2,1)
\]
sends $r\in\mG$ to the infinitesimal generator of the hyperbolic transformations on the hyperboloid $\Hyp\subset\R^{2,1}$ with $r$ as oriented axis,
\item the set $\mG(\gamma)\subset\mG$ is the subset containing the leaves of $\supp(\lambda)$, oriented consistently with the $\lambda$-earthquake whose lifting $\tl$ on $\mH$ fixes $z$, that meet the geodesic arc $[z, h (\gamma)(z)]$,
\item $\mrd\lambda$ denotes $\mrd\meas_\lambda$.
\end{itemize}
Given a smooth family $( h_t)_{t\in I}\subset \mTb$, where $I$ is an interval of $\R$ containing 0, we want to show that for every $\gamma\in\pi_1(\S)$ the map $t\mapsto e^\lambda_l( h_t)(\gamma)$ is smooth. Consider the relative covers $\mH_t$ and subsets $\mG_t(\gamma)\subset\mG$. Denote by $\tl_t$ the realization of $\tl$ in $\mH_t$. Now
\[
e^\lambda_l( h_t)(\gamma)=\int_{\mG}v(r)\chi_{\mG_t(\gamma)}(r)\,\mrd\lambda_t.
\]
For every $t\in I$ there exists a homeomorphism $\zeta_t\colon \de\mH_0\to\de\mH_t$ which is $ h_t$-equivariant, i.e. 
\[
\zeta_t( h_t(\beta)(x))= h_t(\beta)(\zeta_t(x))\ \ \forall x\in\de\mH_0\ \forall\beta\in\pi_1(\S),
\]
and such that for every $x$ that is an endpoint of an axis of $ h_0(\alpha$) for some $\alpha\in\pi_1(\S)$ the map $t\mapsto\zeta_t(x)$ is smooth. It induces a map 
\[
Z_t=\zts\colon \mG\xrightarrow{\sim}\mG.
\]
It turns out that $\lambda_t=Z_t(\lambda_0)$, in the obvious sense. Notice that the endpoints of the leaves of $\lambda_t$ are also endpoints of boundary components for every $t\in I$. Also, $\mG_t(\gamma)(Z_t(s))=\mG_0(\gamma)(s)$ for every $s\in\mG$. Now we have
\begin{align*}
e^\lambda_l( h_t)(\gamma)=&\int_{\mG}v(r)\chi_{\mG_t(\gamma)}(r)\,\mrd Z_t(\lambda_0)=\int_{\mG}v(Z_t(s))\chi_{\mG_0(\gamma)}(s)\,\mrd\lambda_0.
\end{align*}
The integrand of the latter member is a smooth function of $t$, so we get that $t\mapsto e^\lambda_l( h_t)(\gamma)$ is smooth for every $\gamma\in\pi_1(\S)$.
\end{proof}
\end{e_liscia}
{
\tsr
\newtheorem{shu}[r1]{Remark}
\begin{shu}
From the proof of the previous proposition we also get that if $\lambda_n$ is a sequence of laminations converging to $\lambda$ in $\mML$ then $e^{\lambda_n}_l$ converges to $e^\lambda_l$ in $\Gamma(T\mT)$ with the $\mCi$ topology.
\end{shu}
}


\section{Length map}
\label{Sez2}

This section is devoted to find a Hamiltonian $-\bL=-\BL$, given any $\bsl\in\mMLu$, for the vector field $\e{\bsl}=\e{\lambda_1}+\ldots+\e{\lambda_N}$ with respect to a symplectic form on $T\mTb$, provided in the first subsection. After an heuristic computation of $\mrd\bL$ (Subsect. \ref{TsgoBL}), we decompose $\bsl$ in simpler couples still lying in $\mMLu$ such that the sum of their infinitesimal earthquakes gives $\e{\bsl}$ (Subsect. \ref{S_circlam}). For such couples we define $\bL$ (Subsect. \ref{S_Lcirclam}) and show that $\mrd\bL$ is what we expect (Subsect. \ref{Tfov}). Finally, $\BL$ will be constructed as the sum of such length maps (Subsect. \ref{BLgeneric}).

\subsection{The symplectic structure of $\mTb$}
\label{S_sympl}

Fix $\mbb=(b_1,\ldots,b_\n)\in(\R_{>0})^\n$ once for all and consider
\begin{align*}
\mTb=\{h\in\mT\,|\, \ell_h(\de_k)=b_k\ \forall k=1,\ldots,\n\}.
\end{align*}
{A pant decomposition of $\S$ with (internal) curves $\kappa_i$ induces the coordinates 
\[
(\mathbf{l},\boldsymbol\tau,\boldsymbol\beta)=(l_1\ldots,l_{3\g-3+\n},\tau_1,\ldots,\tau_{3\g-3+\n},\beta_1,\ldots,\beta_\n)
\]
on $\mT$, where $l_j$ denotes the length of $\kappa_j$, $\tau_j$ the twist factor of $\kappa_j$, and $\beta_i$ the length of the boundary component $\de_i$ of $\S$. The space $\mTb$ is the submanifold of $\mT$ individuated by the $\n$ equations $\boldsymbol\beta=\mbb$.\\
If $\mu$ has not compact support then there exists $i\in\{1,\ldots,\n\}$ such that $m_i=m(\de_i,\mu)\neq 0$, so we have 
\[
\ell_{E^{t\mu}_l(h)}(\de_i)=|b_i-tm_i|\neq b_i
\]
for $t\in(0,\varepsilon)$ with $\varepsilon$ sufficiently small; such a linear behaviour shows that if $h\in\mTb$ then $e^\mu_l(h)$ does not lie in $T_h\mTb$. However, for every $\bsl$ in
\[
\mMLu=\bigg\{\boldsymbol\mu=(\mu_1,\ldots,\mu_N)\in\mML^{\ N}\,\colon\,\sum_{n=1}^Nm(\de_i,\mu_n)=0,\ i=1,\ldots,\n\bigg\}
\]
$\e{\bsl}$ is a tangent vector field of $\mTb$, as shown at the end of Subsection\ref{S_quakes}.}\\
We need to provide $\mTb$ with a symplectic form $\varpi$. However, there is a natural Weil-Petersson form on $\mTb$ obtained in the following way. Let $2\S$ be the double of $\S$ along its boundary. {Choose a pant decomposition $\kappa_1^\pm,\ldots,\kappa_{6(\g-1)+2\n}^\pm,\de_1,\ldots,\de_\n$ on $2\S$ invariant by the natural involution. Let $\varpi_{{\WP}}$ denote the Weil-Petersson form on the Teichm\"uller space $\mcT_{2\S}$ of $2\S$. It can be written as 
\[
\varpi_{{\WP}}=\sum_{j=1}^{6(\g-1)+2\n} (\mrd\ell^+_j\wedge \mrd\tau^+_j+ \mrd\ell^-_j\wedge \mrd\tau^-_j)+\sum_{i=1}^\n \mrd\ell^0_i\wedge \mrd\tau^0_i
\]
where $\ell_j^\pm$ and $\tau_j^\pm$ denote respectively the length coordinate and the twist coordinate relative to $\kappa_j^\pm$ in $2S$, while $\ell^0_i$ and $\tau^0_i$ denote respectively the length and twist coordinate relative to $\de_i$. Consider the natural immersion $f\colon\mTb\to\mcT_{2\S}$ that doubles a metric on $\S$. With the 2-form
\[
\varpi=f^*\varpi_{{\WP}}=2\sum_{j=1}^{6(\g-1)+2\n} \mrd\ell_j\wedge \mrd\tau_j,
\]
where $\ell_j$ and $\tau_j$ denote respectively the length coordinate and the twist coordinate relative to $f^*(\kappa_j^+)=f^*(\kappa_j^-)$, it turns out that $(\mTb,\varpi)$ is a symplectic manifold.}

\subsection{Hamiltonian conditions}
\label{TsgoBL}

Consider a simple closed curve $\gamma$ not isotopic to a boundary component. Choose a pant decomposition $\{\gamma,\kappa_2,\kappa_3,\ldots\}$ of $\S$. Denoting by $\gamma$ also the measured lamination supported by the curve $\gamma$ with unitary weight, we have for every $h\in\mT$ that
\begin{align*}
\varpi_h(\e{\bsl},e_l^\gamma)&=2\bigg(\mathrm{d}\ell_\gamma\wedge \mathrm{d}\tau_\gamma + \sum_i\mathrm{d}\ell_{\kappa_i}\wedge \mathrm{d}\tau_{\kappa_i}\bigg)(e_l^\gamma,\e{\bsl})=\\
&= \mathrm{d}\ell_\gamma(\e{\bsl})= \mathrm{d}L_\gamma(\e{\bsl})= \sum_{n=1}^N\ddt_{|0}L_\gamma(E^{t\lambda_n}_l(h)).
\end{align*}
Kerckhoff in \cite{K1} proved that on a closed surface $S$ if $\gamma$ and $\nu$ are laminations with a closed curve as support then for every $h$ in the Teichm\"uller space of $S$ the following holds:
\begin{equation}
\label{K1eq}
\ddt_{|0}L_\gamma(E^{t\nu}_l(h))=\int_{\gamma}\cos\theta_{(\gamma,\nu)}(h)\mrd\nu
\end{equation}
where $\theta_{(\gamma,\nu)}(h)$ denotes the angle measured counterclockwise from $\gamma$ to $\nu$ in the $h$-realization. In the proof in \cite{K1} of Equation \eqref{K1eq} the fact that $\nu$ was a closed curve was actually irrelevant. Thus, in our context, the same argument shows that for any $h$ in $\mT$ and $\nu\in\mML$
\[
\ddt_{|0}L_\gamma(E^{t\nu}_l(h))=\int_{\gamma}\cos\theta_{(\gamma,\nu)}(h)\mrd\nu.
\]
Therefore,
\[
\varpi(\e{\bsl},e_l^\gamma)= \sum_{n=1}^N\int_\gamma \cos \theta_{(\gamma,\lambda)} \mrd\lambda_n.
\]
If a function $H\colon \mTb\to\R$ verifies
\[
\mathrm{d}H(e_l^\gamma)=\sum_{n=1}^N\int_\gamma \cos \theta_{(\gamma,\lambda)} \mrd\lambda_n
\]
then, since the space of simple weighted closed curves is dense in $\mCML$, by an approximation argument we get that for every $\nu\in\mCML$ 
\begin{equation}
\label{Hamilt}
\mathrm{d}H(e_l^\nu)=\varpi(\e{\bsl},e_l^\nu).
\end{equation}
Thus, by definition, $H$ is Hamiltonian of the field $\e{\bsl}$.\\
{If $\lambda_1,\ldots,\lambda_N$ have compact support, with the same argument one gets that $H=-\sum_nL_{\lambda_n}$ is a suitable Hamiltonian. In the following sections we will show that it is always possible to construct a Hamiltonian $-\BL$ of $\e{\bsl}$ for every $\bsl\in \mMLu$.}

\subsection{Circuital laminations}
\label{S_circlam}

If $\lambda_1$ and $\lambda_2$ are measured laminations with empty transverse intersection, their sum $\lambda_1\oplus\lambda_2$ is defined by putting $\supp(\lambda_1\oplus\lambda_2)=\supp(\lambda_1)\cup\supp(\lambda_2)$ and $\meas_{\lambda_1\oplus\lambda_2}=\meas_{\lambda_1}+\meas_{\lambda_2}$. By example, if $\lambda=(\delta,\omega)$ is a weighted curve and $\omega=\omega_1+\omega_2$ then $\lambda$ is the sum of $\lambda_1=(\delta,\omega_1)$ and $\lambda_2=(\delta,\omega_2)$. It is immediate to see that 
\begin{equation}
\label{Rcomp}
\e{\lambda_1\oplus\lambda_2}=\e{\lambda_1}+\e{\lambda_2}.
\end{equation}

\tsd
\newtheorem{circuits}[d1]{Definition}
\begin{circuits}
\label{Circuits}
We say that a $I$-uple $\bsm=(\mu_{1},\ldots,\mu_{I})$ of laminations is a circuital lamination if for every $i=1,\ldots,I$
\begin{itemize}
\item $\mu_{i}$ are $\omega$-weighted single spiralling leaves;
\item $\mu_{1},\ldots,\mu_{I}$ are oriented so that for every if $\mu_{i-1}$ ends spiralling near $D_i\in\{\de_1,\ldots,\de_\n\}$ then $\mu_{i}$ starts spiralling near $D_i$, providing $\mu_{0}=\mu_{I}$;
\item the spiralling sense of $\mu_{i-1}$ near $D_i$ is opposite to the one of $\mu_{i}$ near $D_i$.
\end{itemize} 
\end{circuits}

\begin{wrapfigure}{r}{0.5\textwidth}\centering\includegraphics[width=0.9\linewidth]{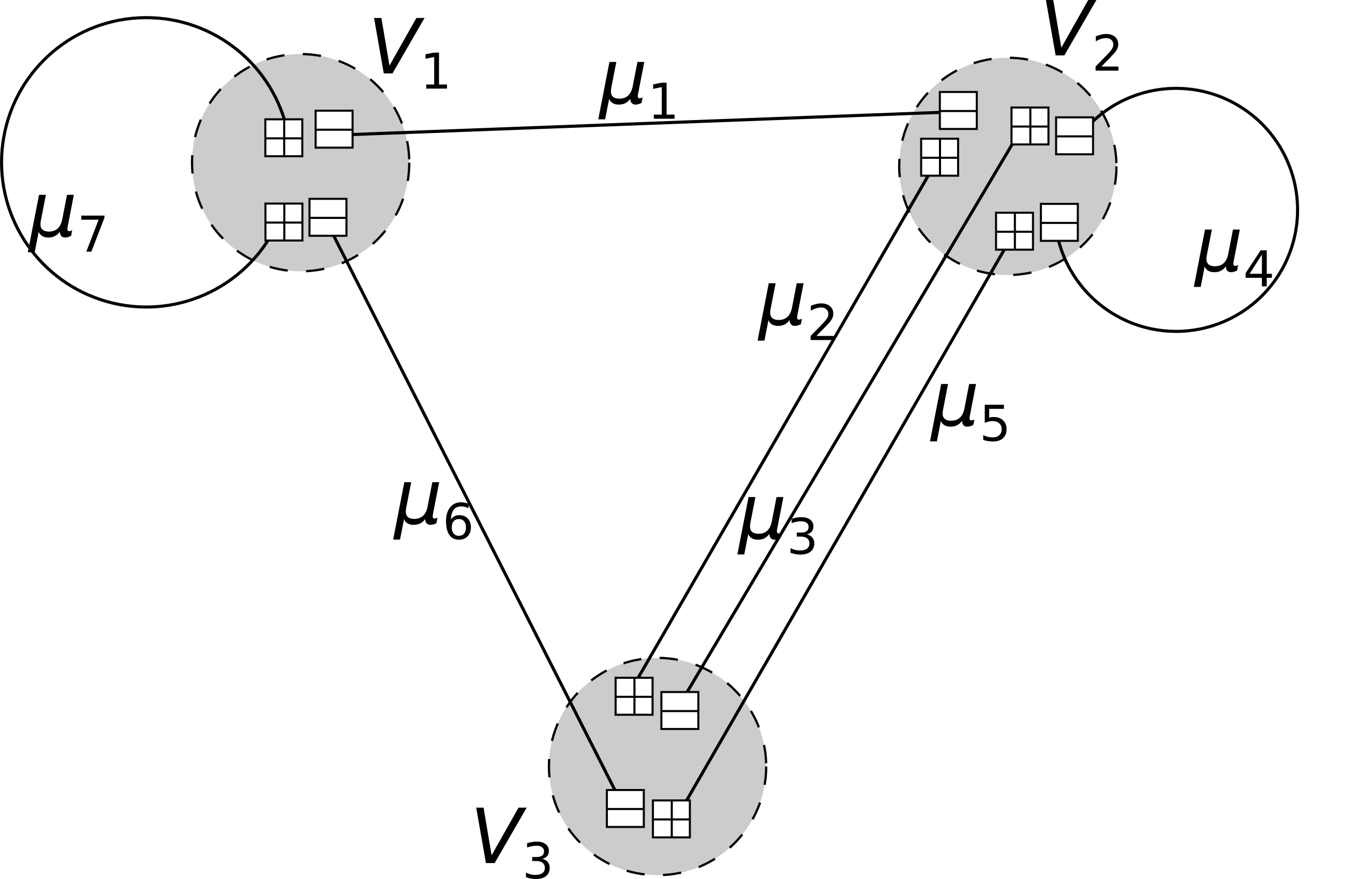}\end{wrapfigure}
A graphic interpretation of such definition can be obtained constructing a multigraph as follows. 
Take $\n$ vertices $V_1,\ldots,V_\n$, representing respectively the boundary components $\de_1,\ldots,\de_\n$ of $S$. For every leaf $\mu_{i}$ spiralling from $\de_m$ to $\de_k$ draw an edge from $V_m$ to $V_k$, marking each endpoint with $\boxminus$ if the leaf spirals in negative sense, with $\boxplus$ otherwise. The $I$-uple $(\mu_{1},\ldots,\mu_{I})$ is circuital if it corresponds to a cycle that every time it passes from an edge to another one switches the sign of the endpoint.

\tsr
\newtheorem{circ_mass}[r1]{Remark}
\begin{circ_mass}
\label{Circ_mass}
If $(\mu_{1},\ldots,\mu_{I})$ is a circuital lamination, then, looking at the corresponding multigraph, for every boundary component $\de_k$ of $\S$
\[
\sum_{i=1}^I(\de_k,\mu_{i})=\omega\cdot \Big(\#\{\boxplus\text{ in }V_k\}-\#\{\boxminus\text{ in }V_k\}\Big)=0.
\]
Therefore, $\mMLu$ contains all the circuital laminations.
\end{circ_mass}

\tsp
\newtheorem{decomp}[p1]{Proposition}
\begin{decomp}
\label{Decomp}
For every $\bsl=(\lambda_1,\ldots,\lambda_N)\in\mMLu$ thre exist circuital laminations $\bsm^{(1)},\ldots,\bsm^{(J)}$ such that 
\begin{equation}
\label{SNBRN}
\e{\bsl}=\e{\bsl^{(0)}}+\sum_{j=1}^J\e{\bsm^{(j)}}
\end{equation}
where $\bsl^{(0)}$ is the $N$-uple of the compact parts of $\lambda_1,\ldots,\lambda_n$.
\begin{proof}
If $\bsl=(\lambda_1^{(0)},\ldots,\lambda_N^{(0)})$ there is nothing to prove. Otherwise, consider the multigraph $G$ associated with $\spir(\bsl)=(\lambda_1\setminus\lambda_1^{(0)},\ldots,\lambda_N\setminus\lambda_N^{(0)})$. We start by looking for a circuital lamination $\bsm^{(1)}=(\mu^{(1)}_{1},\ldots,\mu^{(1)}_{I_1})$ contained in $\spir(\bsl)$; this is equivalent to find a cycle in the graph $G$ alternating the signs of the endpoints of the edges (notice that such cycle is allowed to pass on an edge more than one time).\\
Since $\bsl\in\mMLu$, a vertex $V$ of $G$ contains a $\boxplus$ symbol if and only if $V$ also contains a $\boxminus$ symbol, since the condition $\sum m(\de,\lambda_n)=0$ implies that near $\de$ laminations can not all spiral in the same sense.\\
Let us start from a vertex $D^0$ reached by an endpoint $\boxminus$ of an edge $\hmu_{1}$ and denote by $D^1$ the vertex (maybe coincident with $D^0$) of the other endpoint of $\hmu_{1}$. If such endpoint has the $\boxminus$ symbol, there must be a $\boxplus$ symbol in $D^1$, endpoint of an edge $\hmu_{2}$; vice versa, if such endpoint has the $\boxplus$ symbol, there must be a $\boxminus$ symbol in $D^1$, endpoint of an edge $\hmu_{2}$. Denote by $D_2$ the vertex of the other endpoint of $\hmu_{2}$ and reiterate to find $D^3$ and $\hmu_{3}$, always switching endpoint symbols. Following such procedure, we get a switching path on $G$ (in the sense that consecutive edges have opposite endpoint symbols).\\
If we can find $K$ such that there is $H<K$ and the subpath from $D^H$ to $D^K$ is a switching cycle, then we have finished. We claim that if we visit a vertex $D^k$ for the third time then either we have already found such $K$ (and it is less than $k$) or there is $H<k$ such that the path from $D^H$ to $D^k$ is a switching cycle (so $k$ is the $K$ we were looking for). Suppose we visit a $D^k$ for the third \begin{figure}[h]
\centering
\includegraphics[scale=.10]{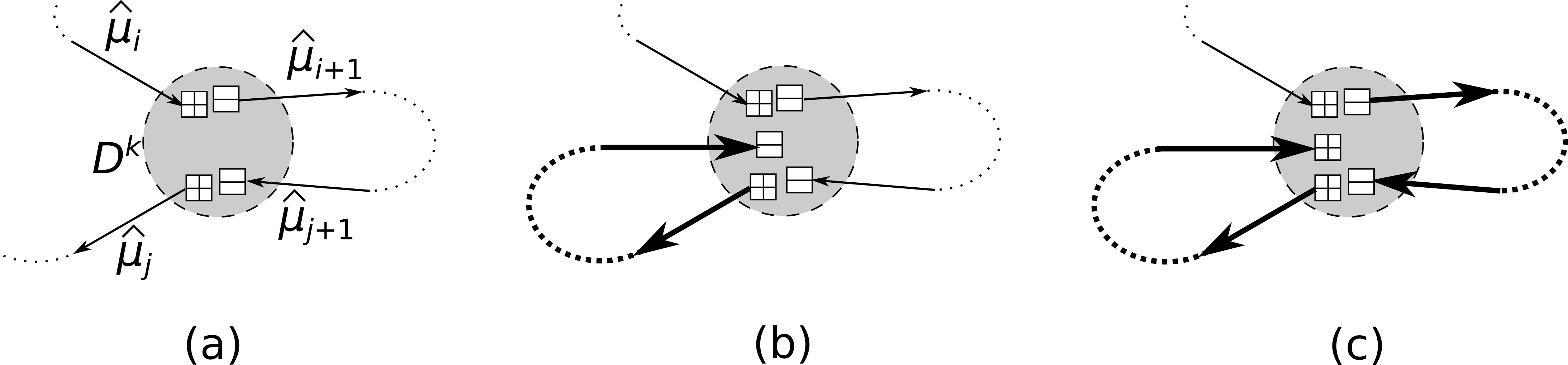}
\caption{}
\label{Dante}
\end{figure}time without having found a switching cycle before. Then the configuration of the previous two visits must be the one in Figure \ref{Dante} (a), up to exchanging $\boxplus$ and $\boxminus$. The third time the path enters $D^k$, it can add either a $\boxminus$ symbol, as in Figure \ref{Dante} (b), or a $\boxplus$ symbol, as in Figure \ref{Dante} (c). In both case, a switching cycle can be individuated, as enlightened in the pictures.\\
So there exists a switching cycle 
\[
(\mu^{(1)}_{1},\mu^{(1)}_{2},\ldots,\mu^{(1)}_{I_1})=(\hmu_{H},\hmu_{H+1},\ldots,\hmu_{K-1})
\]
in $G$, generating a circuital lamination $\bsm^{(1)}$ contained in $\bsl$.\\
We want to endow $\bsm^{(1)}$ with a weight $\omega^{(1)}$ so that if $\bsL=(\Lambda_1,\ldots,\Lambda_N)$ is the $N$-uple of laminations such that 
\[
\e{\bsl}=\e{\bsl^{(0)}}+\e{\bsm^{(1)}}+\e{\bsL}
\]
then $\bsm^{(1)}$ has at least one leaf not contained in the support of $\bsL$. For every spiralling leaf $\delta$ of $\bsl$, denote by $\omega_\delta$ its weight. Define
\[
\omega^{(1)}=\min\bigg\{\frac{\omega_\delta}{\#\{i\in\{1,\ldots,I_1\}\,|\,\mu^{(1)}_{i}=\delta\}}\,\bigg|\,\delta\text{ is a leaf of }\bsl\bigg\}.
\]
In this way, the leaf of $\bsl$ where such minimum is achieved does not appear in the support of $\bsL$.\\
If \begin{figure}[h]
\centering
\includegraphics[scale=.08]{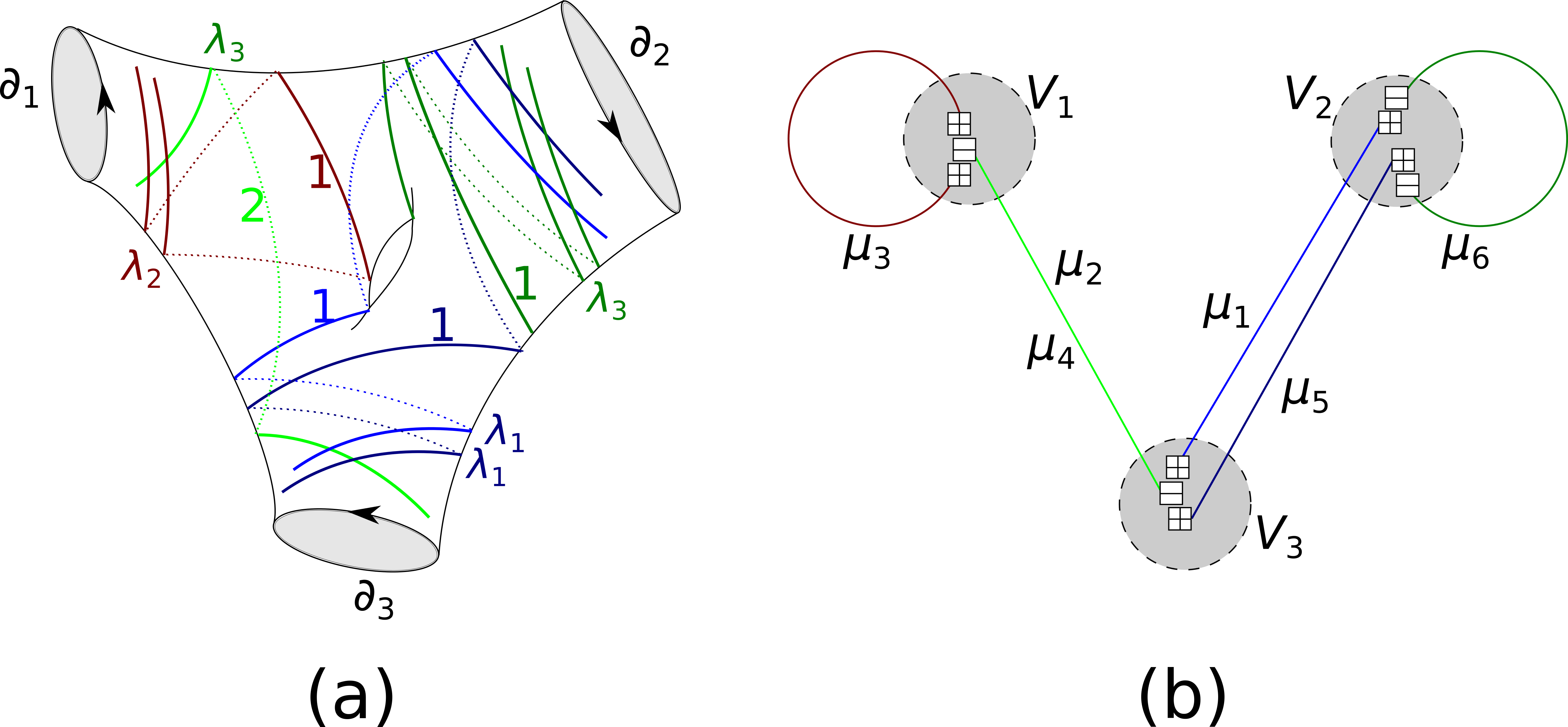}
\caption{}
\label{Circo}
\end{figure}$\bsL$ is the $N$-uple of void laminations, we have finished. See Figure \ref{Circo} as example, where the cycle in (b) spans the triple of laminations in (a). Otherwise, notice that again $\bsL\in\mMLu$ (it depends on the fact that $\bsm^{(1)}$ lies in $\mMLu$; see Remark \ref{Circ_mass}). Moreover $\bsL$ has less leaves than $\bsl$. By a simple inductive argument we get circuital sublaminations $\bsm^{(1)},\ldots,\bsm^{(J)}$, with $J\in\N$, such that \eqref{SNBRN} holds.
\end{proof} 
\end{decomp}

\subsection{The length map for circuital laminations}
\label{S_Lcirclam}

Consider a circuital lamination
\[
\bsl=(\lambda_{1},\ldots,\lambda_{I}).
\]
By definition, there are boundary components $D_1,\ldots,D_I\in\{\de_1,\ldots,\de_\n\}$ of $\S$ such that $\lambda_{i}$ spirals from $D_i$ to $D_{i+1}$ for every $i=1\,\ldots,I$, providing $D_{I+1}=D_{1}$. Also, $\lambda_{i}$ spirals in the opposite sense of $\lambda_{[i+1]}$ near $D_i$.

\tsp
\newtheorem{electro}[p1]{Lemma}
\begin{electro}
\label{De.par}
Let $l$ and $m$ be two simple geodesic in $\S$ spiralling near a boundary component $\de$ in opposite senses, parametrized so that $d(l(t),\de)$ and $d(m(t),\de)$ tend to zero as $t$ goes to infinity. Then there exists a unique $p_0\in l\cap m$ with the following two properties.
\begin{enumerate}
\item \label{Propertyi} Denote by $l_*=l\restr{[t^*,+\infty)}$ and $m_*=m\restr{[T^*,+\infty)}$ the rays in $l$ and $m$ originating at $p_0$ and enumerate consecutively on $l_*$ the elements of $l_*\cap m_*$, starting from $p_0$, as $p_1$, $p_2$, $\ldots$ . Denote by $\hat{l}_k$ the arc of $l_*$ going from $p_k$ to $p_{k+1}$ and by $\hat{m}_k$ the arc of $m_*$ going from $p_k$ to $p_{k+1}$. Then for every $k\in\N$ the piecewise geodesic loop $\hat{l}_k\cup\hat{m}_k$ is isotopic to $\de$.
\item \label{Propertyii} In $l\setminus l_*$ there is no point with the previous property.
\end{enumerate}
\begin{proof}
Clearly, if such $p_0$ exists, then it is unique.\\
On \begin{figure}[h]\centering\includegraphics[scale=.18]{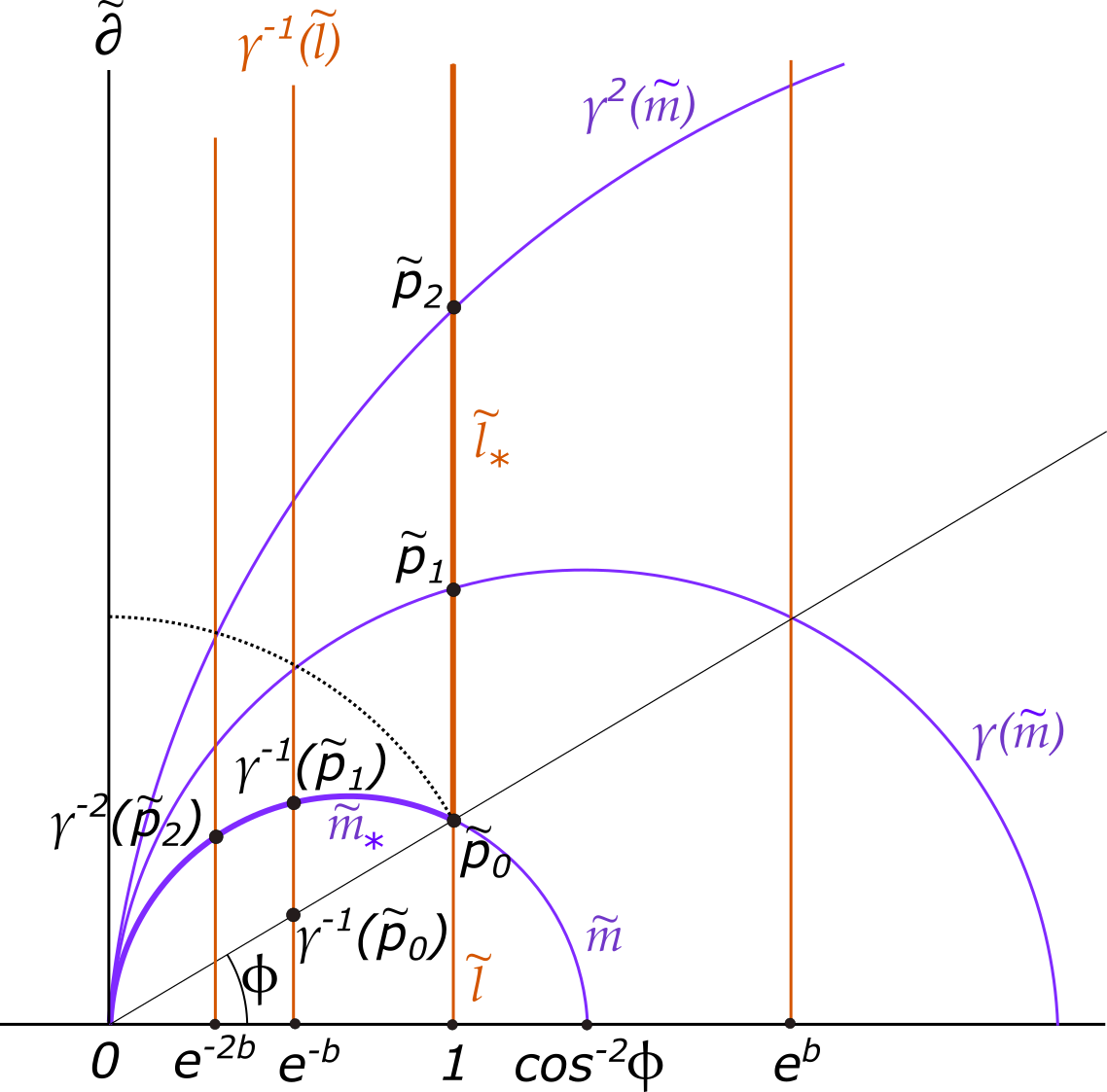}
\end{figure}the universal cover $\mH$ of $\S$ in the upper half-plane model of $\H$ choose coordinates such that a preimage of $\de$ coincides with the imaginary ray and a lift $\tilde{l}$ of $l$ is $1+i\R_{>0}$. Here we are supposing that $l$ spirals around $\de$ in, say, positive sense.\\
Set $b=\ell(\de)$ and let $\gamma\colon z\mapsto e^{b}z$ denote the holonomy transformation corresponding to $\de$. The union of the lifts of $m$ with an ideal endpoint in 0 is $\gamma$-invariant. Among them, there exists a unique $\tilde{m}$ such that $\tilde{l}\cap\gamma^k(\tilde{m})$ is non-empty for every $k\geq 0$ and $\tilde{l}\cap\gamma^k(\tilde{m})$ is empty for every $k<0$. For every $k\geq 0$ let $\tp_k$ be the intersection between $\tilde{l}$ and $\gamma^k(\tilde{m})$ and $p_k$ the projection of $\tp_k$ on $\S$.\\
A simple geometrical analysis shows that $p_0$ satisfies the stated properties.
\end{proof}
\end{electro}

\tsr
\newtheorem{r11}[r1]{Remark}
\begin{r11}
\label{Rtpk}
Let us consider the points $\tp_k$ chosen as in the proof of the previous lemma. They belong to $\tilde{l}$, so $\Re \tp_k=1$ for every $k$. The geodesic $m$ spirals around $\de$ in the opposite sense of $l$, so an ideal endpoint of $\tilde{m}$ must be 0. The other endpoint of $\tilde{m}$ is $\cos^{-2}\phi$, where $\phi=\arg\tp_0$. This implies that $\gamma^k(\tilde{m})$ has ideal endpoints $0$ and $e^{bk}\cos^{-2}\phi$. From this, for every $k\geq 0$ we can compute the imaginary part of the points $\tp_k=l\cap \gamma^k(\tilde{m})$:
\begin{equation*}
\tp_k=\tilde{l}\cap\gamma^k(\tilde{m})=1+i\sqrt{e^{bk}\cos^{-2}\phi-1}.
\end{equation*}
\end{r11}

\tsp
\newtheorem{ppp1}[p1]{Lemma}
\begin{ppp1}
\label{B.collar}
Fix $\mbb\in(\R_{>0})^\n$. For every boundary component $\de$ of $\S$ there exists $\varepsilon(\de)>0$ such that for every $h\in\mTb$ every simple complete geodesic that enters the $\varepsilon(\de)$-collar $\mN(\de)$ of $\de$ exits no more.
\begin{proof}
Choose $h\in\mTb$ and set $b=\ell(\de)$. On the universal cover $\mH\subset \H$ take coordinates such that the imaginary ray projects on a boundary component $\de$. Let $\gamma\colon z\mapsto e^{b}z$ be the corresponding holonomy transformation. {If the endpoints $z<z'$ of a complete geodesic $\tilde{\sigma}$ in $\mH$ are such that $z'>e^bz$, then $z<\gamma(z)<z'<\gamma(z')$, so $\gamma(\tilde{\sigma})$ meets $\tilde{\sigma}$. Therefore, if a geodesic $\sigma\subset\S$ is simple and not spiralling around $\de$, any lift $\tilde{\sigma}$ must have endpoints $z<z'$ such that $z'< e^bz$. A standard computation shows that $\sigma$ does not enter a $\varepsilon(\de)$-collar of $\de$, where 
\[
\varepsilon(\de)=\frac 1{\tanh (b/2)}.
\]
}
\end{proof}
\end{ppp1}

For every boundary component $\de$ of $\S$, we will denote by $\mN(\de)$ the $\varepsilon(\de)$-collar of $\de$ and we will call the union $\mathcal{N}$ of such collars \textit{spiralization neighbourhood}.

\tsr
\newtheorem{rp2}[r1]{Remark}
\begin{rp2}
\label{Rp2}
If $k\geq 1$ then $p_k$ lies in $\mN(\de)$. In fact, a point $x$ of $l$ lies in $\mN(\de)$ if and only if the preimage of $x$ on $\tilde{l}$ has imaginary part greater than $\tan \varphi(\de)=\sinh (b/2)$ (see Lemma \ref{B.collar}). For $k\geq 1$ we have
\[
\Im \tp_k\geq \Im\tp_1=\sqrt{e^{b}\cos^{-2}\phi-1}\geq\sqrt{e^{b}-1}\geq\sinh (b/2).
\]
It may be possible that $p_0$ does not lie in $\mN(\de)$. That is the reason why the definition of $L$ will involve $p_1$ and not $p_0$. 
\end{rp2}

\tsr
\newtheorem{distpk}[r1]{Remark}
\begin{distpk}
\label{Distpk}
{If $k\geq 1$, the distance between $p_k$ and $\de$ is computed by
\[
\tanh d(p_k,\de)=\tanh d(\tp_k, \tde)=\cos\arg \tp_k=\frac{\Re \tp_k}{|\tp_k|}=e^{-bk/2}\cos\phi.
\] }
\end{distpk}

Now let us come back to the circuital lamination $\bsl$ with leaves $\lambda_{1},\ldots,\lambda_{I}$. Let $p^{[i]}=p^{[i]}_{1}$ be the point $p_1$ near $D_i$ chosen as in the proof of Lemma \ref{De.par} when $l=\lambda_{i-1}$ and $m=\lambda_{i}$, providing $\lambda_{0}=\lambda_{I}$. Now we can define a map $L=L_\bsl\colon \mTb\to\R$, that will turn out in Subsection \ref{Tfov} to be the {opposite} of a Hamiltonian of $\e{\bsl}$. 

\tsd
\newtheorem{ldef}[d1]{Definition}
\begin{ldef}
\label{Ldef}
Take an $\omega$-weighted circuital lamination $\bsl$, and consider the points $p^{[i]}$ introduced above. Let $\rho$ be the union of the geodesic arcs in $\lambda_{i}$ with endpoints $p^{[i]}$ and $p^{[i+1]}$ on $i=1,\ldots,I$. For every $h\in\mTb$, set 
\[
L(h)=\omega\bigg\{\ell_h(\rho)+2\log\prod_{i=1}^I\cosh d_h (p^{[i]},D_i)\bigg\}.
\]
\end{ldef}

We notice that $L$ depends on the circuital decomposition of $\bsl$.

\tsr
\newtheorem{lrem}[r1]{Remark}
\begin{lrem}
\label{Lrem}
Consider the loops $\rho_k$ made by the truncations of the leaves $\lambda_{i}$ at the points $p^{[i]}_k$ relative to $D_i$ (defined as in Lemma \ref{De.par}), so that $\rho_1=\rho$. Notice that $\rho_{k+1}\setminus\rho_k$ is a union of $M$ loops, each isotopic to a certain $D_i$. Moreover, such loops tend to some components of $\de\S$, as $k$ goes to infinity. {Setting 
\[
B_h=\sum_{i=1}^I\ell_{h}(D_i),
\]
it turns out that the map
\[
h\mapsto\omega(\ell_h(\rho_k)+2\log\prod_{i=1}^I\cosh d_h(p_k^{[i]},D_i)-kB_h)
\]
is independent on $k$. See \cite{me} for details. Therefore, the map $L_k\colon \mT\to\R$ defined by
\[
L_k(h)=\omega\ell_h(\rho_k)+2\omega\log\prod_{i=1}^I\cosh d_h(p_k^{[i]},D_i)
\]
differs from $L=L_1$ by $(k-1)B_h$, a constant depending only on the $h$-lengths of the boundary components.
}
\end{lrem}

\subsection{The first order variation of $L$}
\label{Tfov}

\begin{figure}[h]
\centering
\includegraphics[scale=.12]{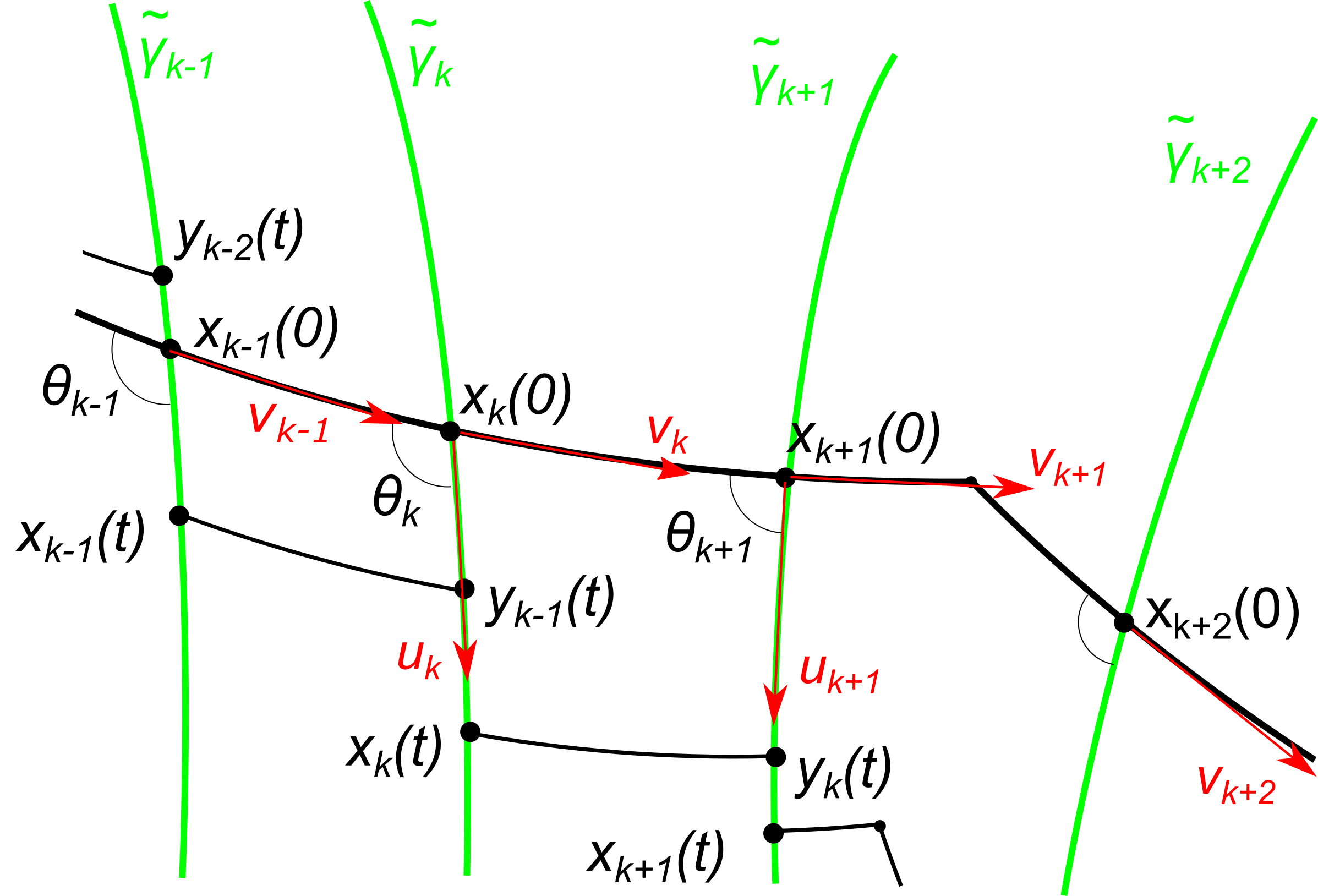}
\caption{}
\label{UkVk}
\end{figure}

The goal of this Subsection is to prove the following proposition:
{
\tsp
\newtheorem{goal.sub}[p1]{Proposition}
\begin{goal.sub}
\label{Goal.sub}
Take an $\omega$-weighted circuit of laminations $\bsl\in\mMLu$ and consider the map $L=L_{\bsl}\colon \mTb\to\R$ given by Definition \ref{Ldef}. For every non-peripheral and non-trivial simple close curve $\gamma$ on $\S$ and for every $h\in\mTb$ the equation
\begin{equation}
\label{Fov.L}
\ddt L(E^{t\gamma}_l(h))=\sum_{i=1}^I\int \cos \theta_{(\lambda_i,\gamma)}(t) \,\mrd\gamma\otimes\mrd\lambda_i
\end{equation}
holds, where $\theta_{(\lambda_i,\gamma)}(t)$ is the angle measured counterclockwise from the support of $\lambda_n$ to $\gamma$, in the $E^{t\gamma}_l(h)$-realization of $\gamma$ and $\lambda_i$.
\end{goal.sub}
Notice that we are slightly abusing the notation, denoting by $\gamma$ also the measured lamination supported by the curve $\gamma$ with unitary weight. This proposition will be true more in general, replacing $\gamma$ with a measured lamination $\nu$ with compact support, as shown at the end of the Subsection.
}\\
Since 
\[
L(h)=\omega\bigg\{\ell_h(\rho)+2\log\prod_{i=1}^I\cosh d_h (p^{[i]},D_i)\bigg\},
\]
we will first compute the derivative in $t=0$ of $\omega\ell_{\Egr(h)}(\rho)$, which will turn out to be
\begin{align*}
\omega\ddt_{|t=0}\ell_{\Egr(h)}(\rho)=\sum_{i=1}^I\int \cos \theta_{(\lambda_i,\gamma)}(0) \,\mrd\gamma\otimes\mrd\lambda_i+\omega\sum_{i=1}^I\mR_i(0)
\end{align*}
where $\mR_i$ are terms due to the presence of the vertices $p^{[i]}$ in $\rho$.\\
After that, setting $F(d)=2\log\cosh d$, we will show that 
\begin{equation}
\label{Cool1}
\mR_i(0)+\ddt_{|t=0}F\big(d_{\Egr(h)}(p^{[i]},D_i)\big)=0
\end{equation}
thus proving Equation \eqref{Fov.L}.\\
\\
Let us start to compute the derivative of $\ell_{\Egr(h)}(\rho)$. Notice that the loop $\rho$ is piecewise geodesic and has exactly $I$ vertices, which are $p^{[i]}$ for $i=1,\ldots,I$.\\
If\begin{figure}[h]
\centering
\includegraphics[scale=.1]{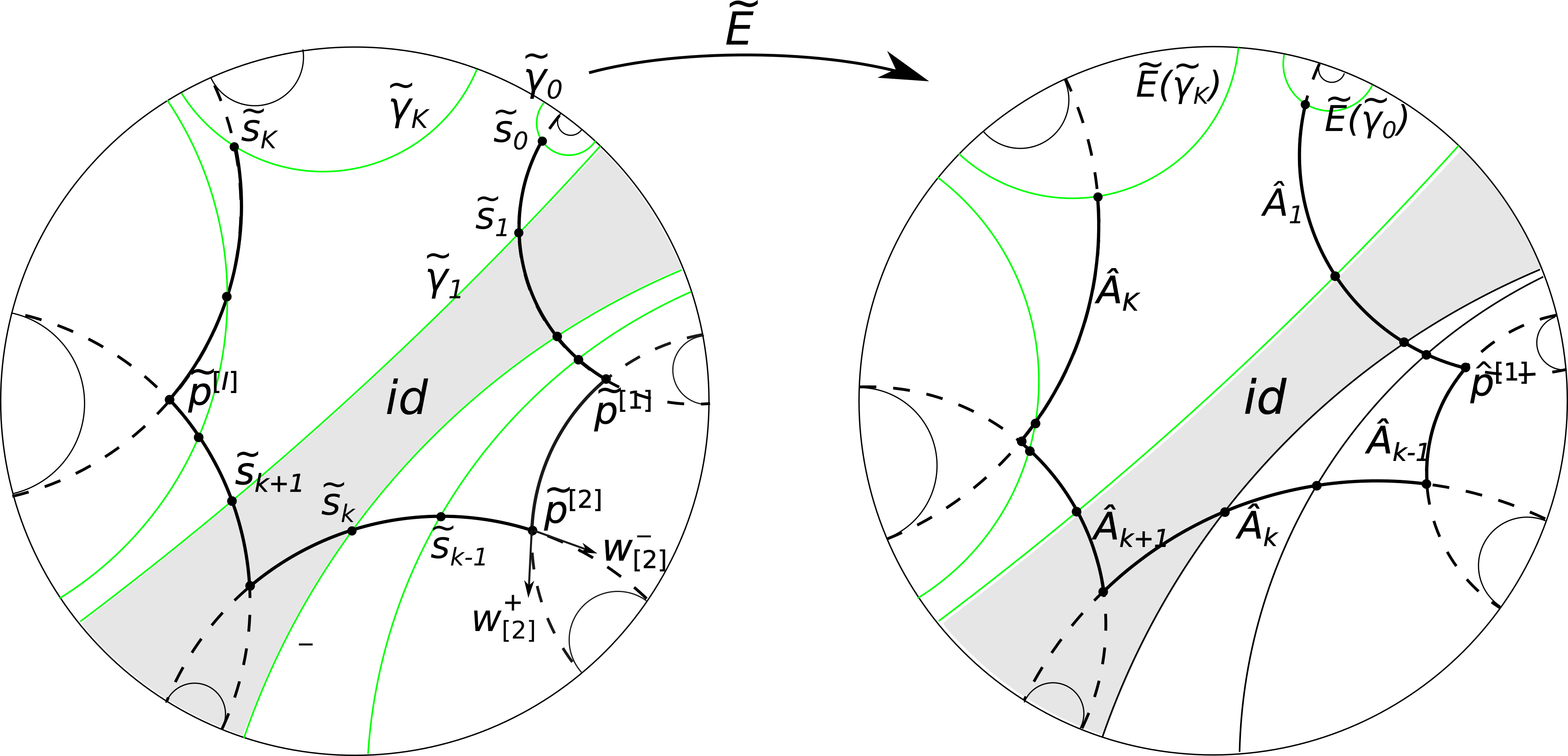}
\caption{Determination of $\hat{\rho}$ and $\hat{A}_j$ (here $I=4$)}
\end{figure} $\iota(\gamma,\lambda_i)=0$ for every $i$ then $\ell(\rho)$ is constant. Otherwise, $\gamma$ meets at least one $\lambda_i$. Notice that $\gamma\cap\rho=\gamma\cap\bigcup\lambda_i$, since $p^{[i]}$ lies in the spiralization neighbourhood for every $i$ (see Lemma \ref{B.collar} and Remark \ref{Rp2}). \begin{figure}[h]
\centering
\includegraphics[scale=.1]{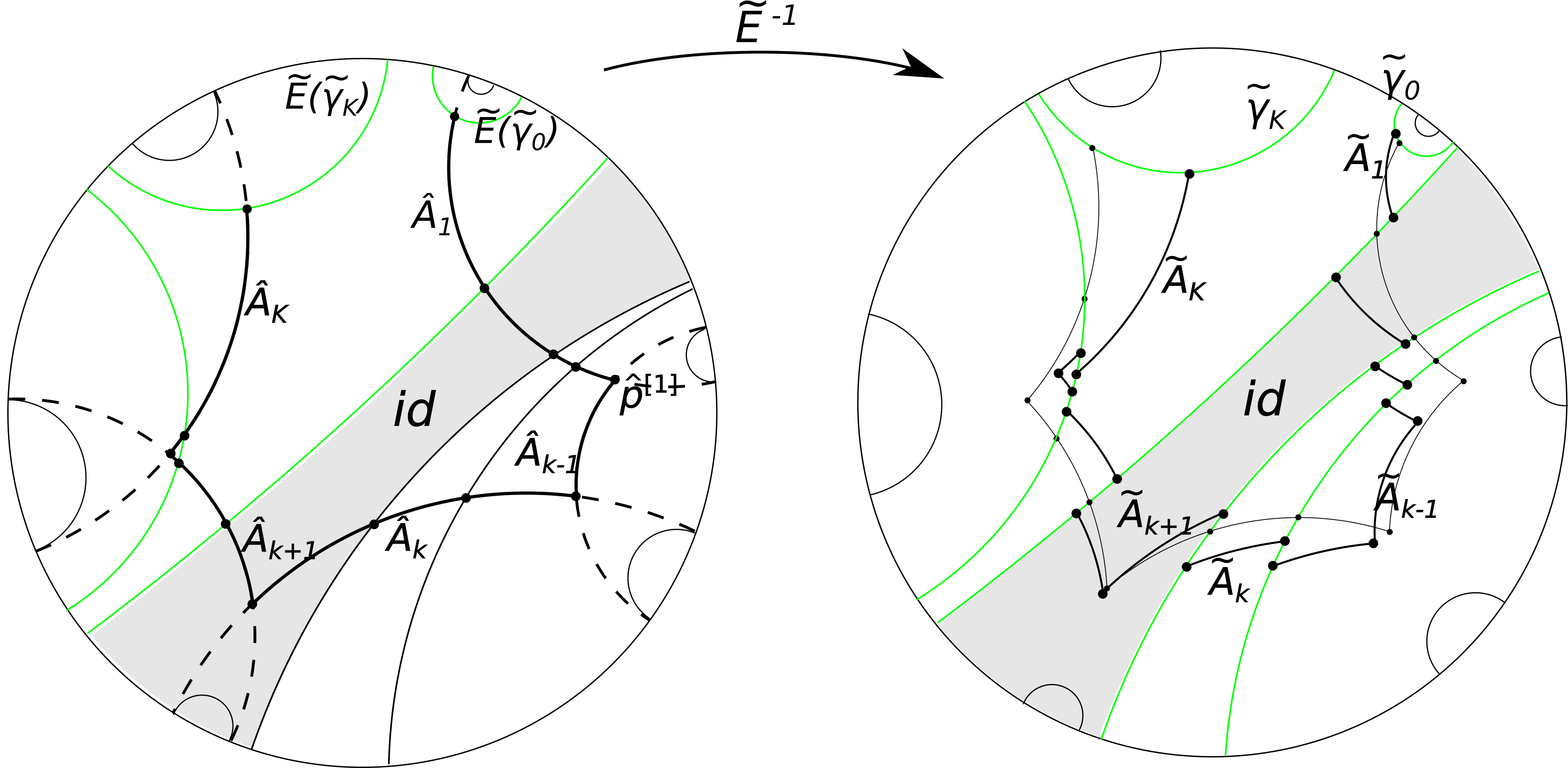}
\caption{Determination of $A_j$ (here $I=4$)}
\end{figure}Choosing an orientation of $\rho$, enumerate consecutively its intersections with $\gamma$ as $s_0,s_1,\ldots,s_{K-1}$. Pick a preimage $\tilde{s}_0$ of $s_0$ on the universal cover $\mH$ of $\S$. If $r\colon [0,1]\to\S$ is a parametrization of the loop $\rho$ such that $r(0)=r(1)=s_0$, take the lift $\tilde{r}\colon [0,1]\to\mH$ with $\tilde{r}(0)=\tilde{s}_0$. Put $\tilde s_K=\tilde r(1)$ and $\tilde s_k$ the preimage of $s_k$ along $\tilde r$ for $k<K$. The preimages of $\gamma$ determine the strata of the lifting $\tE$ of $\Egr$. In particular, denote by $\tg_k$ the preimage of $\gamma$ passing through $\tilde{s}_k$, for $k=0,\ldots,K$.\\
{The path $\tilde{r}$ is piecewise geodesic, with vertices $\tp^{[i]}$. The images of the lifts of the components of $\de\S$ through $\tE$, together with $\tE(\tg_0)$ and $\tE(\tg_K)$, determine the piecewise geodesic arc $\hat{\rho}$ (which does not coincide with $\tE(\tilde{r})$) whose length is equal to $\ell_{\Egr(h)}(\rho)$. The arc $\hat{\rho}$ is divided in $K$ piecewise geodesic subarcs $\hat{A}_1,\ldots,\hat{A}_{K}$ by its intersections with $\bigcup{\tE(\tg_k)}$; such subarcs are enumerated following the orientation of $\hat{\rho}$. The preimage $A_k$ under $\tE$ of $\hat{A}_k$ is a piecewise geodesic arc with endpoints $x_{k}\in\tg_{k-1}$ and $y_k\in\tg_{k}$ with the same length as $\hat{A}_k$. Notice that $x_1=\ts_0$ and $y_{K}=\ts_K$.} This leads to
\[
\ell_{\Egr(h)}(\rho)=\sum_{k=1}^{K}\ell_h(A_k(t)).
\]
For $k=1,\ldots,K$ denote with $v_k$ the unitary vector tangent to $\tilde{r}$ at $\ts_{k}=x_{k-1}(0)=y_{k}(0)$, by $\theta_k$ the angle in $\ts_{k}$ measured counterclockwise from $\tilde{r}$ to $\tg_{k}$ and by $u_k$ the unitary tangent vector to $\tg_{k}$ at $\ts_{k}$ such that $\pi-\theta_k$ is the angle between $v_k$ and $u_k$, as in Figure \ref{UkVk}. Notice that
\[
\sum_{i=1}^I\int\cos\theta_{(\lambda_i,\gamma)}\,\mrd\gamma\otimes\mrd\lambda_i=\omega\sum_{k=1}^K\cos\theta_k.
\]
\tsp
\newtheorem{jb}[p1]{Lemma}
\begin{jb}
\label{Jb}
For $k=1,\ldots,K-1$, the following identity holds:
\begin{equation}
\label{Altern}
\dot{x}_{k+1}(0)=\dot{y}_{k}(0)+u_k.
\end{equation}
\begin{proof}
Denote by $d_k(t)$ the signed distance between $y_{k}(0)=x_{k+1}(0)$ and $y_{k}(t)$ on $\tg_{k}$ oriented as $u_k$. Then
\begin{gather*}
y_{k}(t)=y_k(0)\cosh d_k(t)+u_k\sinh d_k(t)\\
x_{k+1}(t)=x_{k+1}(0)\cosh \big(d_k(t)+t\big)+u_k\sinh \big(d_k(t)+t\big).
\end{gather*}
Therefore,
\begin{gather*}
\dot{y}_{k}(0)=u_k\dot{d}_k(0)\\
\dot{x}_{k+1}(0)=u_k\big(\dot{d}_k(0)+1\big)
\end{gather*}
leading to \eqref{Altern}.
\end{proof} 
\end{jb}
\begin{figure}[h]
\centering
\includegraphics[scale=.1]{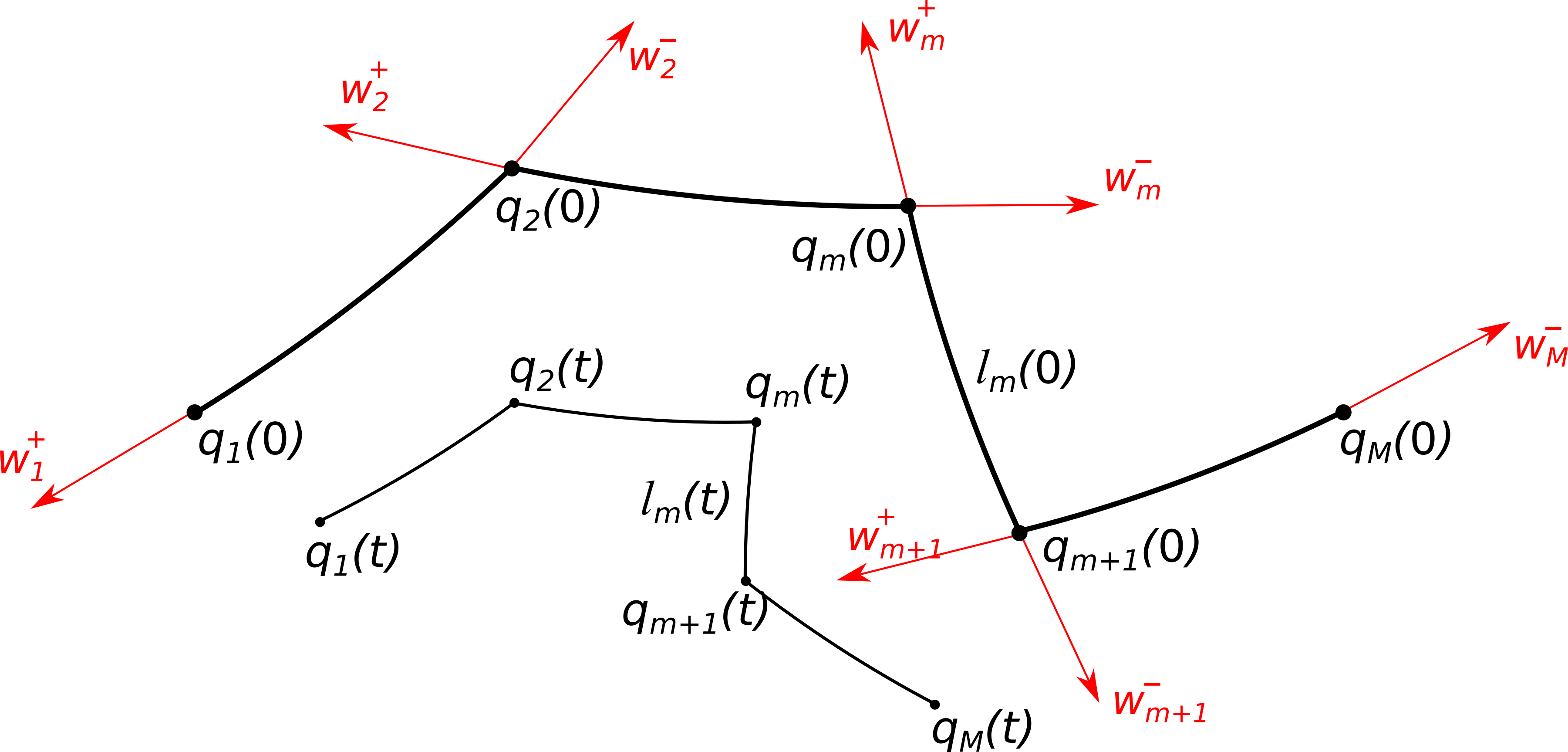}
\caption{}
\label{Polychain}
\end{figure}
\tsp
\newtheorem{chain_fov}[p1]{Lemma}
\begin{chain_fov}
\label{Chain_fov}
Consider the hyperboloid model of $\H$ in $\R^{2,1}=(\R^3,\langle *,*\rangle)$ (where $\langle x,y\rangle=-x_0y_0+x_1y_1+x_2y_2$), namely 
\[
\H\cong\{x\in\R^{2,1}\colon \langle x,x\rangle=-1,x_0>0\}.
\]
Given an integer $M\geq 2$ and a $\mC^1$ map $q\colon [0,1]\to(\H)^M$, let $C(t)$ be the oriented open polygonal chain in $\H$ of vertices $q_1(t)\ldots,q_M(t)$. Denote by $w_{m}^-$ and $-w_{m}^+$ respectively the left and right unitary tangent vector to $C(0)$ at $q_m(0)$, for $m=2,\ldots,M-1$. Define analogously $-w^+_1$ and $w^-_M$, as in Figure \ref{Polychain}. Then 
\[
\ddt_{|0}\ell(C(t))=\langle\dot{q}_1(0),w_1^+\rangle+\langle\dot{q}_M(0),w_M^+\rangle+\sum_{m=2}^{M-1}\langle\dot{q}_m(0),w_m^-+w_m^+\rangle.
\]
\begin{proof}
Set $l_m(t)=d(q_m(t),q_{m+1}(t))$. It suffices to prove that 
\begin{equation}
\label{Pnc}
\dot{l}_m(0)=\langle\dot{q}_m(0),w_m^+\rangle+\langle\dot{q}_{m+1}(0),w_{m+1}^-\rangle
\end{equation}
for $m=1,\ldots,M-1$.
Since $\cosh l_m(t)=-\langle q_m(t),q_{m+1}(t)\rangle$, differentiating at $t=0$ we get
\begin{equation}
\label{Delight}
\dot{l}_m(0)\sinh l_m(0)=-\langle \dot{q}_m(0),q_{m+1}(0)\rangle-\langle q_{m}(0),\dot{q}_{m+1}(0)\rangle.
\end{equation}
Since
\begin{gather*}
q_{m+1}(0)=q_{m}(0)\cosh l_j(0)-w_m^+\sinh l_m(0)\\
q_{m}(0)=q_{m+1}(0)\cosh l_j(0)-w_{m+1}^-\sinh l_m(0),
\end{gather*}
equation \eqref{Delight} becomes
\[
\dot{l}_m(0)\sinh l_m(0)=\langle\dot{q}_m(0),w_m^+\rangle\sinh l_m(0)+\langle\dot{q}_{m+1}(0),w_{m+1}^-\rangle\sinh l_m(0),
\]
which gives \eqref{Pnc}.
\end{proof}
\end{chain_fov}
We are able now to prove the following result.
\tsp
\newtheorem{sistem1}[p1]{Proposition}
\begin{sistem1}
\label{Sistem1}
\[
\ddt_{|0}\sum_{k=1}^K\ell_h(A_k(t))=\sum_{k=1}^{K}\cos\theta_k+\sum_{i=1}^I\mR_i(0)
\]
where $\mR_1,\ldots,\mR_I$ are terms related to the $I$ vertices of $\rho$ (explicitly computed in the proof, see Equation \eqref{S1}).
\begin{proof}
Each $A_k$ is a piecewise geodesic arc, with endpoints $x_k$ and $y_k$. Applying Lemma \ref{Chain_fov} to every $A_k$, we get 
\begin{align*}
\sum_{k=1}^K\ddt_{|0}&\ell_h(A_k(t))=\sum_{k=1}^K\Big(\langle \dot{x}_k(0),-v_{k-1}\rangle+\langle \dot{y}_k(0),v_{k}\rangle\Big)+\\
&+\sum_{i=1}^I\big\langle \dot{\tp}^{[i]}(0),w_-^{[i]})+w_+^{[i]}\big\rangle,
\end{align*}
where $\pm w_\pm^{[i]}$ denote the unitary vectors tangent to $\tilde{r}$ at $\tp^{[i]}$ and the vectors $v_k$ where defined before Lemma \ref{Jb}.\\
Let us put
\begin{gather}
\mS=\sum_{k=1}^K\Big(\langle \dot{x}_k(0),-v_{k-1}\rangle+\langle \dot{y}_k(0),v_{k}\rangle\Big)\nonumber\\
\label{S1}
\mR_i(0)=\big\langle \dot{\tp}^{[i]}(0),w_-^{[i]}+w_+^{[i]}\big\rangle
\end{gather}
Using \eqref{Altern}, we have that
\begin{align*}
\mS=&\langle \dot{x}_1(0),-v_0\rangle+\sum_{k=2}^K\langle \dot{x}_k(0),-v_{k-1}\rangle+\sum_{k=1}^{K-1}\langle \dot{y}_k(0),v_{k}\rangle+\langle \dot{y}_K(0),v_{K}\rangle=\\
=&-\langle \dot{x}_1(0),v_0\rangle+\sum_{k=1}^{K-1}\langle \dot{x}_{k+1}(0)-\dot{y}_k(0),-v_{k}\rangle+\langle \dot{y}_K(0),v_{K}\rangle=\\
=&-\langle \dot{x}_1(0),v_0\rangle+\sum_{k=1}^{K-1}\langle u_k,-v_{k}\rangle+\langle \dot{y}_K(0),v_{K}\rangle.
\end{align*}
Since $s_0$ is a point were $\rho$ is smooth and $\tilde{s}_0=x_1(0)$ and $\tilde{s}_{m}=y_{K}(0)$ are preimages of $s_0$, there exists a covering transformation $T$ such that $T\dot{x}_1(0)=\dot{y}_K(0)+u_K$ and $Tv_0=v_K$. Now
\begin{align*}
\mS=&-\langle T\dot{x}_1(0),Tv_0\rangle-\sum_{k=1}^{K-1}\langle u_k,v_{k}\rangle+\langle \dot{y}_K(0),v_{K}\rangle=\\
=&-\langle \dot{y}_K(0)+u_K,v_K\rangle-\sum_{k=1}^{K-1}\langle u_k,v_{k}\rangle+\langle \dot{y}_K(0),v_{K}\rangle=\\
=&-\sum_{k=1}^{K}\langle u_k,v_{k}\rangle=\sum_{k=1}^{K}\cos\theta_k.
\end{align*}
\end{proof}
\end{sistem1}

Now we have to show that Equation \eqref{Cool1} holds.\\
{Let us first recall some known facts on the hyperboloid model of $\Hyp$, keeping the notation of the proof of Proposition \ref{Sistem1}. For every geodesic $\gamma$ in $\Hyp$ there is a space-like vector $n_\gamma$ such that
\[
\gamma=\{\ulx\in\Hyp\,|\,\langle \ulx,n_\gamma\rangle=0\}.
\]
The boundary at infinity of $\Hyp$ is identified with 
\[
\dinf\Hyp=\quot{\{\ulx\in\R^{2,1}\,|\,\langle\ulx,\ulx\rangle=0\}}{ \ulx\sim a\ulx, a\in\R^*}
\]
and its elements will be written within square brackets. See also \cite{B-P}.\\
There is a notion of cross product in $\R^{2,1}$, analogous to the Euclidean environment: if $\mrd V$ denotes the volume form in $\R^{2,1}$, the cross product between $\ulx\in\R^{2,1}$ and $\uly\in\R^{2,1}$ is the vector $\ulx\boxtimes\uly\in\R^{2,1}$ such that for every $\ulz\in\R^{2,1}$
\[
\langle \ulx\boxtimes\uly,\ulz\rangle=\mrd V(\ulx,\uly,\ulz).
\]
The following hold:
\begin{gather*}
\langle\ulx,\uly\boxtimes\ulz\rangle=\langle\ulz,\ulx\boxtimes\uly\rangle\\
(\ulx\boxtimes\uly)\boxtimes\ulz=\langle\uly,\ulz\rangle\ulx-\langle\ulx,\ulz\rangle\uly\\
\langle \ulx\boxtimes\uly,\ulx\boxtimes\uly\rangle=\langle\ulx,\uly \rangle^2-\langle\ulx,\ulx \rangle\langle\uly,\uly \rangle
\end{gather*}
for every $\ulx,\uly,\ulz\in\R^{2,1}$.}\\
\begin{wrapfigure}{r}{0.3\textwidth}\centering\includegraphics[width=0.9\linewidth]{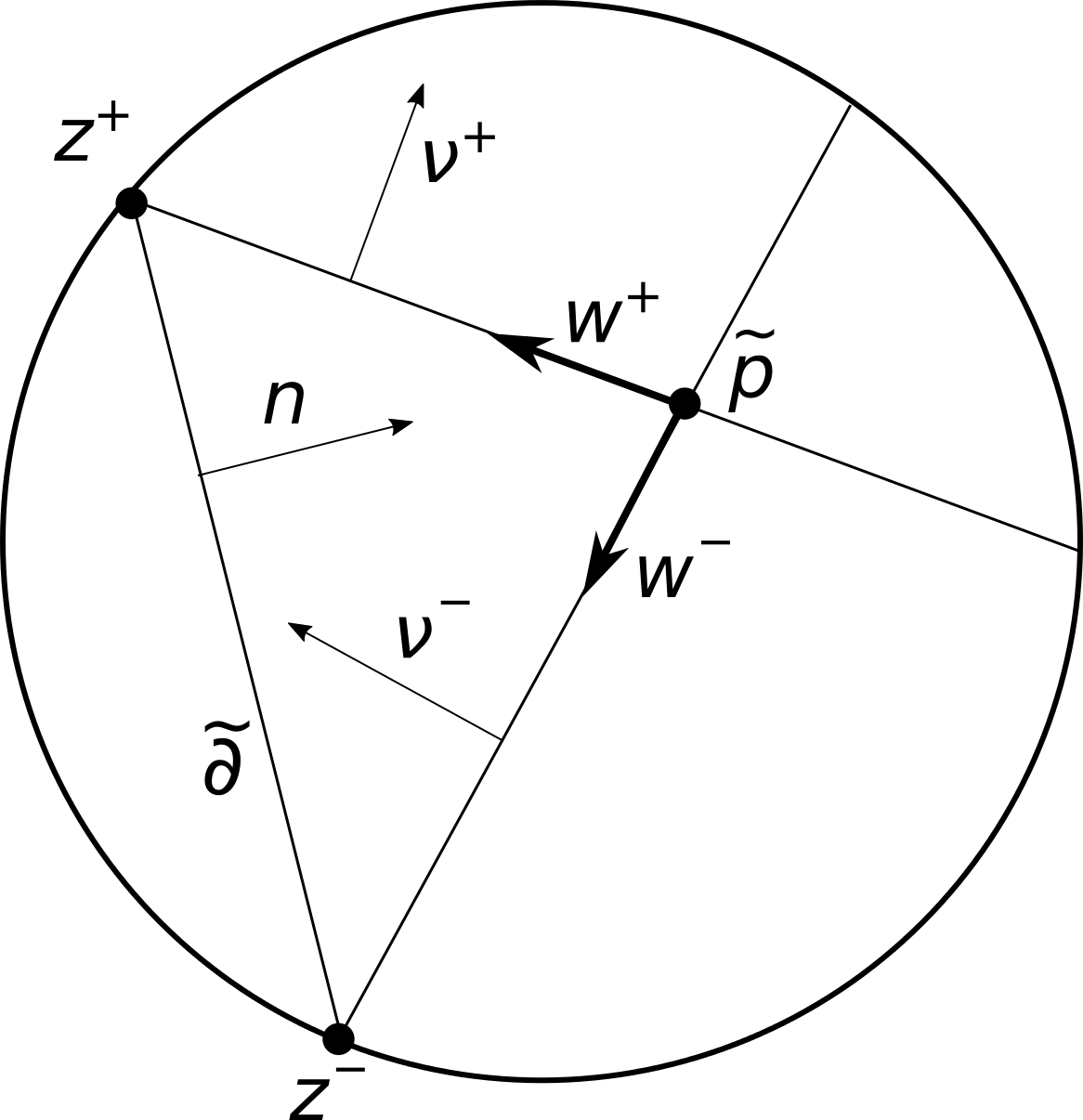}\end{wrapfigure}
Now fix $i$ and denote by $\de$ the component of $\de\S$ whose spiralization neighbourhood contains $p^{[i]}$. If $\tde$ is the lift of $\de$ closer to $\tp=\tp^{[i]}$, denote by $[z^+]$ and $[z^-]$ the ideal endpoints of $\tde$, so that $w^\pm=w_\pm^{[i]}$ is pointing towards $[z^\pm]$.\\ 
The unitary vector
\[
n=\frac{z^-\boxtimes z^+}{\|z^-\boxtimes z^+\|_{2,1}}
\]
is the normal unitary vector of $\tde$ pointing towards $\tp$. %
Up to precomposing by a proper isometry, we can suppose that $[z^+]$ and $[z^-]$ are kept fixed by $\tE$, thus $\tE(n)=n$. If $d=d(p,\de)=d(\tp,\tde)$, then $\sinh d=\langle \tp,n\rangle$. Therefore
\[
\dot{d}=\frac{\langle\dot{\tp},n\rangle}{\cosh d}
\]
and
\[
\ddt F(d)=2\frac{\sinh d}{\cosh^2 d}\langle\dot{\tp},n\rangle,
\]
where we have set $F(d)= 2\log\cosh d$. Now Equation \eqref{Cool1} becomes
\[
\bigg\langle\dot{\tp},w^++w^-+2\frac{\sinh d}{\cosh^2 d}n\bigg\rangle=0.
\]
The following proposition will prove such equation computing $w^\pm$ in terms of $\tp$ and $n$.

\tsp
\newtheorem{pcool}[p1]{Proposition}
\begin{pcool}
\[
\langle\dot{\tp},w^++w^-\rangle=-2\frac{\sinh d}{\cosh^2 d}\langle\dot{\tp},n\rangle.
\]
\begin{proof}
The vector $w^\pm$ can be written as $\tp\boxtimes \nu^\pm$, where $\nu^\pm$ is the unitary vector tangent to $\H$ and normal to $w^\pm$ (i.e. to $\lambda_{[i-1]}/\lambda_{[i]}$) oriented in the proper way; namely,
\[
\nu^\pm=-\frac{z^\pm\boxtimes \tp}{\|z^\pm\boxtimes \tp\|_{2,1}}.
\]
Thus,
\[
w^\pm=-\tp\boxtimes\frac{z^\pm\boxtimes \tp}{\|z^\pm\boxtimes \tp\|_{2,1}}=-\frac{-\langle \tp,\tp\rangle z^\pm+\langle z^\pm,\tp\rangle \tp}{\langle z^\pm, \tp\rangle}=-\frac{z^\pm+\langle z^\pm,\tp\rangle \tp}{\langle z^\pm, \tp\rangle}.
\]
{We claim that 
\begin{equation}
\label{zulp}
z^\pm=\tp-(\sinh d)n\pm \tp\boxtimes n.
\end{equation}
First, we have to see that the right hand side of \eqref{zulp} is a null vector; let us compute the square norm of $\tp\boxtimes n$:
\[
\langle\tp\boxtimes n,\tp\boxtimes n\rangle=\langle\tp, n\rangle^2-\langle \tp,\tp\rangle \langle n,n\rangle=\sinh^2d+1=\cosh^2d
\]
Now
\begin{align*}
&\langle \tp-(\sinh d)n\pm \tp\boxtimes n,\tp-(\sinh d)n\pm \tp\boxtimes n\rangle=\\
=&\langle\tp,\tp \rangle -(\sinh d)\langle \tp,n\rangle -(\sinh d)\langle n,\tp\rangle +(\sinh^2 d)\langle n,n\rangle+\langle \tp\boxtimes n,\tp\boxtimes n\rangle=\\
=&-1-\sinh^2d-\sinh^2d+\sinh^2d+\cosh^2d=0.
\end{align*} 
On the other hand, we have to check that $\tp-(\sinh d)n\pm \tp\boxtimes n$ are the ideal endpoints of $\tde$ (or equivalently $\langle \tp-(\sinh d)n\pm \tp\boxtimes n,n\rangle=0$) such that 
\[
\Big(\langle\tp-(\sinh d)n- \tp\boxtimes n,\ \tp-(\sinh d)n+ \tp\boxtimes n,\  n\Big)
\]
forms a negative basis of $\R^{2,1}$. Now
\[
\langle\tp-(\sinh d)n\pm \tp\boxtimes n,n\rangle=\langle\tp,n\rangle-(\sinh d)\langle n,n\rangle=0
\]
and
\begin{align*}
&\langle\tp-(\sinh d)n- \tp\boxtimes n,(\tp-(\sinh d)n+ \tp\boxtimes n)\boxtimes n\rangle=\\
=&\langle \tp-(\sinh d)n-\tp\boxtimes n,\tp\boxtimes n + \tp\rangle=-1-\cosh^2d<0.
\end{align*}
}
Thus, we can compute
\[
\langle z^\pm, \tp\rangle=\langle \tp-(\sinh d)n\pm \tp\boxtimes n, \tp\rangle=-\cosh^2d
\]
and
\[
w^\pm=-\frac{z^\pm+\langle z^\pm,\tp\rangle \tp}{\langle z^\pm,\ul p\rangle}=-\frac{(\sinh^2d)\tp+(\sinh d)n\mp \tp\boxtimes n}{\cosh^2d}.
\]
Now
\begin{align*}
\langle \dot{\tp},w^++w^-\rangle=\bigg\langle\dot{\tp},-\frac{(2\sinh^2d)\tp+(2\sinh d)n}{\cosh^2d}\bigg\rangle=-2\frac{\sinh d}{\cosh^2d}\langle \dot{\tp},n\rangle.
\end{align*}
\end{proof}
\end{pcool}

{
Finally, let us consider the first order variation of $t\mapsto L(E^{t\nu}_l(h))$ in the general case, when $\nu\in\mCML$.
\tsp
\newtheorem{gen.fov}[p1]{{Proposition}}
\begin{gen.fov}
Consider a circuital lamination $\bsl\in\mMLu$. For every $h\in\mTb$ and $\nu\in\mCML$ the following formula holds:
\begin{equation}
\label{Tfovgen}
\ddt L_{\bsl}(E^{t\nu}_l(h))=\sum_{i=1}^I\int \cos \theta_{(\lambda_i,\nu)}(t) \mrd\nu\otimes\mrd\lambda_i.
\end{equation}
\begin{proof}
The space of weighted curves on $\S$ is dense in $\mCML$ (see \cite{P-H}), so take a sequence $(\gamma_j)$ of weighted curves converging to $\nu$. With the notation
\[
\Cos(\lambda,\mu)(t)=\int \cos \theta_{(\lambda,\mu)}(t) \mrd\mu\otimes\mrd\lambda
\]
used in \cite{K2}, we have seen that
\[
\ddt L_{\bsl}(E^{t\gamma_j}_l(h))=\sum_{i=1}^I\Cos(\lambda_i,\gamma_j)(t).
\]
Clearly $L_{\bsl}(E^{t\gamma_j}_l(h))_{|t=0}=L_{\bsl}(E^{t\nu}_l(h))_{|t=0}$ for every $j$, so if we prove that $\sum_i\Cos(\lambda_i,\gamma_j)$ tends uniformly to $\sum_i\Cos(\lambda_i,\nu)$ then $L_{\bsl}(E^{t\gamma_j}_l(h))$ tends to $L_{\bsl}(E^{t\nu}_l(h))$ and \eqref{Tfovgen} holds. Kerckhoff showed in \cite{K2} itself that $\Cos(\delta,\gamma_j)$ tends uniformly to $\Cos(\delta,\nu)$ for every $\delta$ closed curve in $\S$, {but his argument still works if $\delta$ is a spiralling leaf of a lamination on $\S$}, so we can conclude.  
\end{proof}
\end{gen.fov}
}

\subsection{The map $\BL\colon \mTb\to\R$}
\label{BLgeneric}

In order to extend the definition of $L_\bsl$ to any $\bsl=(\lambda_1,\ldots,\lambda_N)\in\mMLu$, consider the decomposition 
\[
\e{\bsl}=\e{\bsl^{(0)}}+\sum_{j=1}^J\e{\bsm^{(j)}}
\]
by Proposition \ref{Decomp}, where $\bsl^{(0)}$ is the compact part of $\bsl$ and $\bsm^{(j)}$ are circuital laminations. Now define 
\begin{equation}
\label{Eris}
\BL=\sum_{j=1}^JL^{(j)}
\end{equation}
where $L^{(0)}=\sum_nL_{\lambda_n^{(0)}}$ and $L^{(j)}$ is the length map of $\bsm^{(j)}$ in Definition \ref{Ldef}, for $j\neq 0$. Since Equation \eqref{Fov.L} holds for every $L^{(j)}$, we can deduce
\begin{equation*}
\ddt\BL(E^{t\nu}_l(h))=\sum_{n=1}^N\int\cos\theta_{(\lambda_n,\nu)}\mrd\nu\otimes\mrd\lambda
\end{equation*}
for every $\nu\in\mCML$ and $h\in\mTb$. In particular, $-\BL$ is a Hamiltonian of $\e{\bsl}$ (see Equation \eqref{Hamilt}).


\section{Properties of $\BL$}
\label{Sez3}

\subsection{$\BL$ is proper}
\label{BLprop}

{Now we are going to show that the map $\BL$ is proper under the hypothesis that $\bsl=(\lambda_1,\ldots,\lambda_N)$ \textit{fills up} $\S$, which means that every non-trivial non-peripheral simple closed curve on $\S$ meets $\bigcup\supp(\lambda_n)$. Set
\[
\mFMLu=\{\bsl\in\mMLu\,|\, \text{$\bsl$ fills up $S$}\}\cup\{\vv\}.
\]
As explained in section \ref{S_topML} any spiralling geodesic $\gamma$ of a measured geodesic lamination can be replaced by a geodesic arc $\gamma^R$ orthogonal to the boundary. For each $\nu\in\mML$ denote by $\nu^R$ the set of geodesic arcs obtained by $\nu$ replacing each spiralling geodesic $\gamma$ of $\nu$ with $\gamma^R$ and set $\bsl^R=\big((\lambda_1)^R,\ldots,(\lambda_N)^R\big)$. Notice that if $\bsl\in\mFMLu\setminus\{\vv\}$ then $\bsl^R$ still fills up $\S$.

\tsp
\newtheorem{filling.ter}[p1]{Lemma}
\begin{filling.ter}
\label{Filling.ter}
Consider two disjoint geodesics $\de$ and $\de'$ in $\H$, a geodesic $\gamma$ going from an endpoint of $\de$ to an endpoint of $\de'$, the geodesic arc $\gamma^R$ with endpoints on $\de$ and $\de'$ normal to $\de$ and $\de'$, two positive real numbers $\epsilon,\epsilon'\leq\ell(\gamma^R)/2$, the $\epsilon$-collars $N$ of $\de$ and the $\epsilon'$-collar $N'$ of $\de'$. Then
\[
\ell(\gamma\setminus(N\cup N'))\geq\ell(\gamma^R\setminus (N\cup N'))=\ell(\gamma^R)-\epsilon-\epsilon'.
\]
\qed
\end{filling.ter}

Such lemma is quite easy to prove; see \cite{me} for details.

\tsp
\newtheorem{BL.propria}[p1]{Proposition}
\begin{BL.propria}
\label{BL.Propria}
If $\bsl\in\mFMLu\setminus\{\vv\}$ then the map $\BL\colon \mTb\to\R$ is proper.
\begin{proof}
Choose a pant decomposition of $\S$ with curves $\kappa_1,\ldots,\kappa_{3(\g-1)+\n}$, $\de_1,\ldots,\de_\n$ and consider the related coordinates 
\[
(l_1,\ldots,l_{3(\g-1)+\n},\tau_1,\ldots,\tau_{3(\g-1)+\n})
\]
on $\mTb$, where $l_i$ is the length of $\kappa_i$ and $\tau_i$ is the twist factor on $\kappa_i$. Choose also for every $\kappa_i$ two dual curves $\kappa^*_i$ and $\kappa_i^{**}$ whose lengths can reconstruct $\tau_i$ (as explained in \cite{F-L-P}; see Figure \ref{Dualcurves}).\\
We have seen at the beginning of this subsection that if $\bsl\in\mFMLu$ then $\bsl^R$ fills up $\S$; \begin{figure}[h]
\centering
\includegraphics[scale=.10]{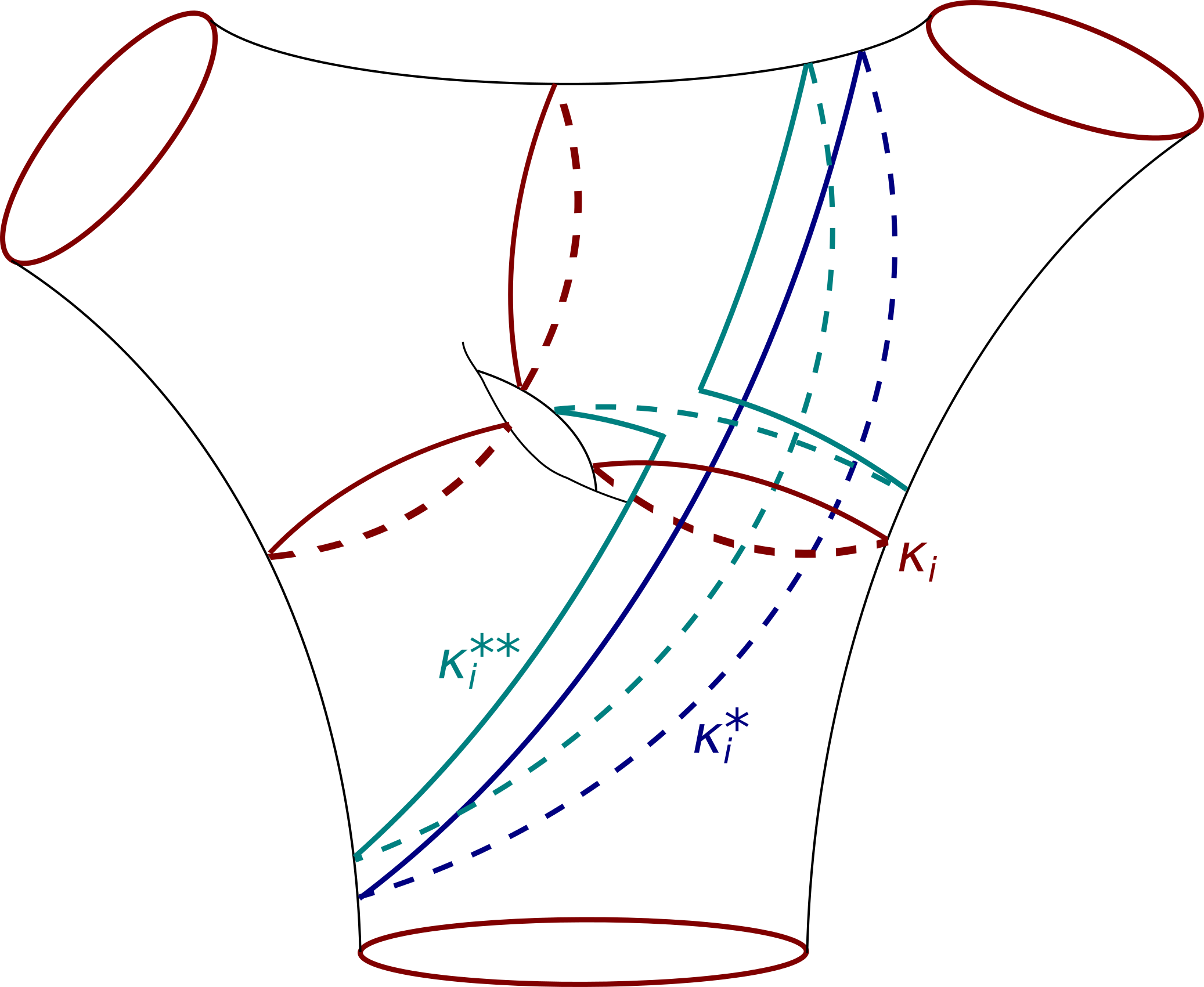}
\caption{}
\label{Dualcurves}
\end{figure}this implies 
that every simple closed non-trivial curve in $\S$ is isotopic to a curve on $G=D\cup\bigcup\supp\big((\lambda_n)^R\big)$, where $D=\bigcup \de_j$.\\
We claim that 
\begin{align*}
\bL_G\colon &\mTb\longrightarrow\R\\
&h\mapsto \sum_{j=1}^\n b_j+\sum_{n=1}^N\ell_h\big((\lambda_n)^R\big)
\end{align*}
is a proper map. Pick a divergent sequence $\{h_k\}$ in $\mTb$; then the sequence
\[
\{(l_1,\ldots,l_{3(\g-1)+\n},\tau_1,\ldots,\tau_{3(\g-1)+\n})(h_k)\}
\]
is divergent in $\R^{6(\g-1)+2\n}$. This implies that
\[
S_k=\sum_{i=1}^{3(\g-1)+\n}\big[\ell_{h_k}([\kappa_i])+\ell_{h_k}([\kappa^*_i])+\ell_{h_k}([\kappa^{**}_i])\big]\xrightarrow{n\to\infty} +\infty,
\]
where for any closed curve $\kappa$ and hyperbolic metric $h$ we denote by $\ell_h([\kappa])$ the $h$-length of the geodesic $h$-realization of $\kappa$.\\
{Each $\kappa_i$ (and $\kappa^*_i$ and $\kappa^{**}_i$) is isotopic to many (not necessarily simple) curves in $G$, but for every $i$ the number
\begin{align*}
m_i=\min\bigg\{&\max_{p\in G}\Big\{\#\big(\pi^{-1}(p)\cap([0,1]\times\{0\})\big)\Big\}\,\Big|\,\pi\colon [0,1]\times[0,1]\to \S \text{ isotopy} \\
&\text{between $\pi(*,0)=\kappa_i$ and $\pi(*,1)$ closed curve in }G\bigg\},
\end{align*}
which denotes a sort of minimum of the degrees of the isotopies between $\kappa_i$ and any curve in $G$, does not depend on the metric.} The same holds for $m^*_i$ and $m^{**}_i$ (the analogous numbers for $\kappa^*_i$ and $\kappa^{**}_i$ respectively). If $m_0$ is the maximum among all $m_i$'s, $m^*_i$'s and $m^{**}_i$'s, then 
\[
S_k\leq 3m_0(3\g-3+\n)\bL_G(h_k).
\]
Therefore, $\{\bL_G(h_k)\}$ is going to infinity as $\{h_k\}$ is diverging.\\
Since $\bL_G(h_k)=\sum b_i+\sum\ell_{h_k}\big((\lambda_n)^R\big)$ is diverging, two possibilities occur:
\begin{itemize}
\item a compact sublamination $\boldsymbol\gamma^R$ of $\bsl^R$ has divergent length; but since $\boldsymbol\gamma^R=\boldsymbol\gamma$, also $\BL(h_k)$ is diverging;
\item no compact sublamination of $\bsl^R$ has divergent length; then an arc $\gamma^R$ in $\bsl^R$ (replacement of a spiralling leaf $\gamma$ of $\bsl$ between $\de$ and $\de'$) has divergent length. Also $\ell_{h_k}\big(\gamma\setminus (\mN(\de)\cup \mN(\de'))\big)$ diverges, by Lemma \ref{Filling.ter}, where $\mN(\de)$ is the $\varepsilon(\de)$-collar introduced in Subsection \ref{S_Lcirclam}. 
From the definition, 
\[
\BL(h_k)>\omega \ell_{h_k}(\gamma-\mN(\de)-\mN(\de'))>\omega\Big(\ell_{h_k}(\gamma^R)-\varepsilon(\de)-\varepsilon(\de')\Big),
\]
implying that $\BL(h_k)$ is diverging. 
\end{itemize}
\end{proof}
\end{BL.propria}

\subsection{The second order variation of $\BL$}
\label{SovBL}

If $\bsl\in\mFMLu\setminus\{\vv\}$, the map $\BL\colon\mTb\to\R_{\geq 0}$ is strictly convex along left earthquakes, which means that $t\mapsto \bL(E^{t\nu}_l(h))$ is strictly convex for every $\nu$ in $\mCML$ and every $h\in\mTb$. Kerckhoff proved it in \cite{K1} for $\nu\in\mCML$, but his argument still applies to spiralling laminations: he worked on the universal cover of the surface $\tilde{S}$, where the key-point was that any right (respectively left) earthquake induces a homeomorphism on $\dinf\tilde{S}$ that moves clockwise (respectively counterclockwise), which still is true in our context.

\tsr
\newtheorem{min_point}[r1]{Remark}
\begin{min_point}
\label{Min_point}
As in \cite{K1}, properness and strict convexity of $\bL\colon\mTb\to\R$ assure that $\bL$ admits exactly one point of minimum $h_0$.
\end{min_point}

The goal of this subsection is to show that the Hessian of $\BL$ is positive definite on a critical point $h_0\in\mTb$ of $\BL$. If $\lambda_1,\ldots,\lambda_N$ have compact discrete support, then the result is already known through explicit formulas (see \cite{W2}, \cite{Br}), which however involve quantities that are not meaningful in our setting. Let us consider $\nu\in\mCML$. We already know from Subsections \ref{Tfov} and \ref{BLgeneric} that
\[
\ddt\BL(E^{t\nu}_l(h_0))=\sum_{n=1}^N\int \cos \theta_{(\lambda_n,\nu)}(t) \,\mrd\nu\otimes\mrd\lambda_n
\]
holds, where $\theta_{(\lambda_n,\nu)}(t)$ is the angle measured counterclockwise from the support of $\lambda_n$ to $\nu$, in the $E^{t\nu}_l(h_0)$-realization of $\nu$ and $\lambda_n$.\\
The compact part of $\nu$ is approximated by closed weighted curves, so let us consider first a unitary closed curve $\gamma$.
If $\delta$ is a weighted spiralling leaf of $\bsl$, we will first compute
\[
\ddt_{|0}\int \cos \theta_{(\delta,\gamma)}(t) \,\mrd\gamma\otimes\mrd\delta=\sum_{i=1}^m\ddt_{|0}\cos\theta_i,
\]
where, enumerating consecutively along $\delta$ the points $x_1,\ldots,x_m$ in $\zeta\cap\gamma$, $\theta_i$ is the angle measured counterclockwise from $\delta$ to $\gamma$ at $x_i$. Then we will deduce an estimate which guarantees that even passing at the limit of closed curves the second derivative stays positive.\\
\begin{wrapfigure}{r}{0.35\textwidth}\centering\includegraphics[width=0.9\linewidth]{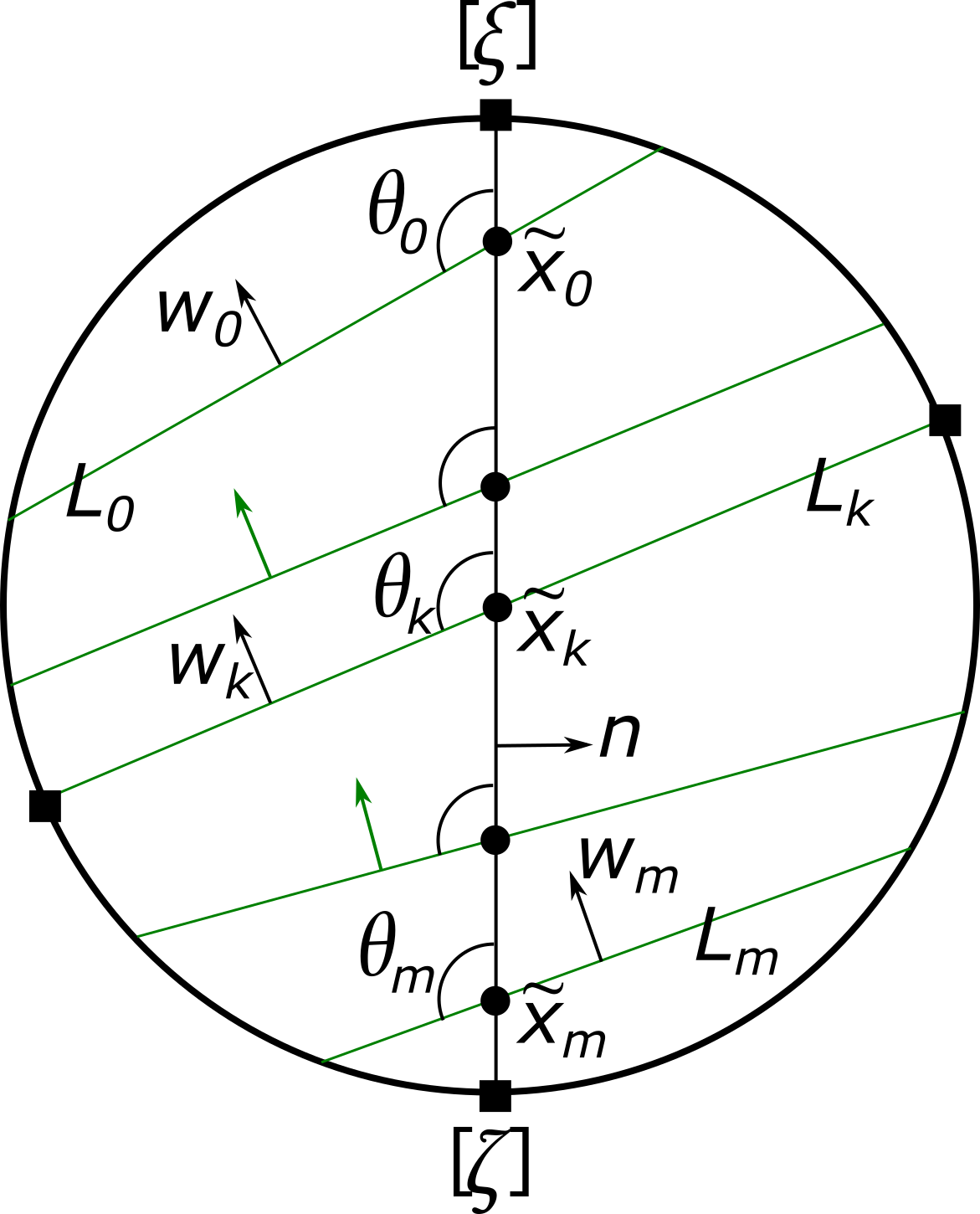}
\end{wrapfigure}
Let us transfer the problem on the universal covering $\mH\subset\Hyp$ of $\S$ in the hyperboloid model of $\H$. Fix a lift $\td$ of $\delta$; denote by $\tx_1,\ldots,\tx_m$ the preimages of $x_1,\ldots,x_m$ on $\td$ and by $L_1,\ldots,L_m$ the liftings of $\gamma$ passing respectively through $\tx_1,\ldots,\tx_m$. Denote by $[\xi]$ and $[\zeta]$ the ideal endpoints of $\td$ so that $\tx_1,\ldots,\tx_m$ are enumerated from $[\xi]$ to $[\zeta]$ and $\xi_0=\zeta_0=1$, if we write vectors $\ulx$ in $\R^{2,1}$ as $\ulx=(x_0,x_1,x_2)$. We can choose coordinates such that $\langle\xi,\zeta\rangle=-1$. Fix $k\in\{1,\ldots,m\}$ and consider the lift $\tE^t$ of $E^{t\gamma}_l$ which fixes the gap whose boundary contains $L_k$ and $L_{k-1}$ (if $k=1$ take the earthquake that fixes the gap adjacent with $L_1$ whose ideal boundary contains $[\xi]$). Choose unitary vectors $w_1,\ldots,w_m$ normal respectively to $L_1,\ldots,L_m$ so that $\cos\theta_i=\langle w_k,n\rangle$ for every $i$. 
Now, since we are in the hyperboloid model of $\H$, let us identify $\R^{2,1}$ with the Lie algebra $\mathfrak{so}(2,1)$. Now
\begin{gather*}
\xi(t)=\tE^t(\xi)=\exp(-tw_1)\cdots\exp(-tw_{k-1})\xi,\\
\zeta(t)=\tE^t(\zeta)=\exp(+tw_k)\cdots\exp(+tw_m)\zeta,\\
n(t)=\frac{\xi(t)\boxtimes\zeta(t)}{\|\xi(t)\boxtimes\zeta(t)\|_{2,1}}=\frac{\xi(t)\boxtimes\zeta(t)}{-\langle\xi(t),\zeta(t)\rangle}
\end{gather*}
so
\begin{gather}
\label{dxi}
\dot{\xi}(0)=-\sum_{i=1}^{k-1}w_i\boxtimes\xi,\\
\label{deta}
\dot{\zeta}(0)=\sum_{i=k}^{m}w_i\boxtimes\zeta.
\end{gather}
Since 
\[
\ddt_{|0}\cos\theta_k(t)=\ddt_{|0}\langle w_k,n(t)\rangle=\langle w_k,\dot{n}(0)\rangle,
\]
let us compute $\dot{n}(0)$. In general,
\begin{align*}
\dot{n}(0)=&\frac{\dot{\xi}(0)\boxtimes\zeta+\xi\boxtimes\dot{\zeta}(0)}{-\langle\xi\boxtimes\zeta\rangle}+\frac{\xi\boxtimes\zeta}{\langle\xi\boxtimes\zeta\rangle^2}\ddt_{|0}\langle\xi(t),\zeta(t)\rangle=\\
=&\dot{\xi}(0)\boxtimes\zeta+\xi\boxtimes\dot{\zeta}(0)+n\cdot\ddt_{|0}\langle\xi(t),\zeta(t)\rangle.
\end{align*}
Setting $z=\dot{\xi}(0)\boxtimes\zeta+{\xi}\boxtimes\dot{\zeta}(0)$, we deduce that there is $\beta\in\R$ such that $\dot{n}(0)=z+\beta n$. So from
\[
0=\langle\dot{n}(0),n\rangle=\langle z,n\rangle +\beta\langle n,n\rangle=\langle z,n\rangle +\beta
\]
we get 
\[
\dot{n}(0)=z-\langle z,n\rangle n.
\]
Writing $\dot{\xi}$ for $\dot{\xi}(0)$ and $\dot{\zeta}$ for $\dot{\zeta}(0)$, setting for every $i$
\begin{equation*}
w_i=a_i\xi+b_i\zeta+c_i n,
\end{equation*}
(notice that $a_i>0>b_i$) and using \eqref{dxi}, \eqref{deta}, we compute $z$ as
\begin{align*}
z=&\dot{\xi}\boxtimes\zeta+\xi\boxtimes\dot{\zeta}=-\sum_{i=1}^{k-1}(w_i\boxtimes\xi)\boxtimes\zeta+\sum_{i=k}^{m}\xi\boxtimes(w_i\boxtimes\zeta)=\\
=&\sum_{i=1}^{k-1}b_i\zeta+\sum_{i=k}^{m}a_i\xi+\sum_{i=1}^m c_i n.
\end{align*}
Now 
\begin{align*}
\ddt_{|0}\cos\theta_k(t)=&\langle w_k,\dot{n}(0)\rangle=\langle w_k,z\rangle-\langle z,n\rangle\langle w_k,n\rangle.
\end{align*}
The three products take values
\begin{align*}
\langle w_k,z\rangle= & \bigg\langle a_k\xi+b_k\zeta+c_k n\ ,\ \sum_{i=1}^{k-1}b_i\zeta+\sum_{i=k}^{m}a_i\xi+\sum_{i=1}^m c_i n\bigg\rangle=\\
=&-\sum_{i=1}^{k-1}a_kb_i-\sum_{i=k}^{m}a_ib_k+\sum_{i=1}^m c_ic_k\\
\langle z,n\rangle=&\bigg\langle\sum_{i=1}^{k-1}b_i\zeta+\sum_{i=k}^{m}a_i\xi+\sum_{i=1}^m c_i n\ ,\ n\bigg\rangle=\sum_{i=1}^mc_i\\
\langle w_k,n\rangle=&\langle a_k\xi+b_k\zeta+c_k n\ ,\ n\rangle=c_k
\end{align*}
so 
\[
\ddt_{|0}\cos\theta_k(t)=-\sum_{i=1}^{k-1}a_kb_i-\sum_{i=k}^{m}a_ib_k.
\]
The sum over $k$ gives
\[
\sum_{k=1}^m\ddt_{|0}\cos\theta_k(t)=-\sum_{k=1}^m\sum_{i=1}^{k-1}a_kb_i-\sum_{k=1}^m\sum_{i=k}^{m}a_ib_k=-\sum_{k=1}^ma_kb_k-2\sum_{i<k}a_ib_k.
\]
Notice that $c_k=\langle w_k,n\rangle=\cos\theta_k$ and 
\[
1=\langle w_k,w_k\rangle = -2a_kb_k+c_k^2,
\]
which implies $-a_kb_k=(\sin^2\theta_k)/2$. The terms $r_{ik}=-a_ib_k>0$ have the property that $r_{ik}r_{ki}=(\sin^2\theta_i\sin^2\theta_k)/4$; moreover,
\[
\cosh d(L_i,L_k)=\langle w_i,w_k\rangle=-a_ib_k-a_kb_i+c_ic_k=r_{ik}+r_{ki}+\cos\theta_i\cos\theta_k.
\]
Since $d(L_i,L_k)$ is bounded by the maximal length of a curve contained in $\bigcup\supp(\lambda_n)\setminus\mN$, there is $M_0>0$ such that
\[
r_{ik}+r_{ki}=\cosh d(L_i,L_k)-\cos\theta_i\cos\theta_k\leq \cosh M_0 +1.
\]
Now
\[
r_{ik}=\frac{r_{ik}r_{ki}}{r_{ki}}\geq \frac{r_{ik}r_{ki}}{r_{ik}+r_{ki}}\geq \frac{\sin^2\theta_i\sin^2\theta_k}{4(\cosh M_0+1)}.
\]
We finally get 
\[
\sum_{k=1}^m\ddt_{|0}\cos\theta_k(t)\geq \frac 12\bigg(\sum_{k=1}^m\sin^2\theta_k+\sum_{i<k}\frac{\sin^2\theta_i\sin^2\theta_k}{2(\cosh M_0+1)}\bigg).
\]
This holds for a spiralling leaf $\delta$ in $\bsl$. Considering all the leaves of $\bsl$, there is $M_1>0$ such that we obtain
\[
\frac{\mrd^2}{\mrd t^2}\bL(E^{t\gamma}_l(h_0))\geq M_1\sum_{n=1}^N\int_{\lambda_n}\int_{\lambda_n}\sin^2\theta_{(\lambda_n,\gamma)}(x)\sin^2\theta_{(\lambda_n,\gamma)}(y)\mrd\gamma(x)\mrd\gamma(y)
\]
or equivalently
\[
\Hess\bL(e^\gamma_l,e^\gamma_l)\restr{h_0}\geq M_1\sum_{n=1}^N\iint\sin^2\theta_{(\lambda_n,\gamma)}(x)\sin^2\theta_{(\lambda_n,\gamma)}(y)\mrd\gamma(x)\mrd\gamma(y).
\]
Now let us consider a generic $\nu\in\mCML$. It is the limit of weighted closed curves $\gamma_k$. As for the first order variation of $\BL$, with an approximation argument we get  
that
\[
\Hess\bL(e^\nu_l,e^\nu_l)\restr{h_0}\geq M_1\sum_{n=1}^N\iint\sin^2\theta_{(\lambda_n,\nu)}(x)\sin^2\theta_{(\lambda_n,\nu)}(y)\mrd\nu(x)\mrd\nu(y).
\]
Therefore, $\Hess_{h_0}\BL$ is definite positive.

\tsr
\newtheorem{emily}[r1]{Remark}
\begin{emily}
\label{Emily}
When $N=2$, consider $\bsl=(\lambda_1,\lambda_2)\in\mFMLu\setminus\{\vv\}$ and let $h_0$ be the unique critical point of $\BL$. Since $-\BL$ is the symplectic gradient of $\e{\bsl}$ with respect to $\varpi$, as shown in Section \ref{Sez2}, saying that $\Hess_{h_0}\BL>0$ is equivalent to state that $e^{\lambda_1}_l$ and $-e^{\lambda_2}_l=e^{\lambda_2}_r$ meet transversely (only) in $h_0$. \\
If $\bsl\in\mMLu\setminus\mFMLu$, with $N\geq 2$, the map $\BL$ does not have any critical points. Kerckhoff proved it for $\n=0$ and $N=2$ in \cite{K2} (Theorem 2.1 II), but the same argument works for any $\n\geq 0$ and $N\geq 2$, since the key point was that Equation \eqref{Tfovgen} holds.
\end{emily}


\section{{The tangent space}}
\label{Sez4}

In this section we extend a result achieved in \cite{B-S}, Appendix B., using the enlightened properties of the Hamiltonian $\BL$ of the vector field $e^\lambda_l+e^\mu_l$.\\ 
Let $\Mink=(\R^{2,1},\langle*,*\rangle)$ be the 3-dimensional Minkowski space and consider $\H$ as the set of unitary future-pointing vectors of $\Mink$, preserved by the action of $\SO$. In this section, we identify $\mT$ with the space of cocompact Fuchsian representations $h\colon\piS\to \SO$ up to conjugacy. An affine deformation of $h$ is a representation 
\[
\rho=h+\tau\colon\piS\to\R^3\rtimes\SO\subset\Isom(\Mink)=\R^3\rtimes O(2,1)
\]
with $\tau(\gamma)\in\R^3$ for every $\gamma\in\piS$. The space $T\mT$ is identified with  
\[
\{\rho=h+\tau\colon\piS\to\R^3\rtimes\SO\,\colon\,h\in\mT\}/\text{conj.},
\]
the space of affine deformations of Fuchsian representations.\\
In \cite{Ba3}, Barbot proved that for every $\rho\in\mRb$ there are two maximal disjoint convex non-empty domains $\Omega^\pm(\rho)\subset\Mink$ such that
\begin{itemize}
\item[-] $\Omega^+(\rho)$ (respectively $\Omega^-(\rho)$) is complete in the future (respectively in the past);
\item[-] the action of $\rho(\piS)$ on $\Omega^\pm(\rho)$ is free and properly discontinuous;
\item[-] $\rho(\piS)\backslash\Omega^\pm(\rho)\simeq S\times\R$ is a maximal Cauchy-hyperbolic spacetime. 
\end{itemize}
Being $\Omega^\pm(\rho)$ regular domains, they are associated with two measured laminations $\lambda_\pm$, considered as dual to the singular loci of $\Omega^\pm(\rho)$. See \cite{B-B} for details. Denote by $\Psi:T\mT\to\mML^{\ 2}$ the arising map. As in \cite{B-S}, if $\Psi\hta=\lpm$ then 
\begin{equation}
\label{saens}
\tau=e^{\pl}_r(h)=e^{\ml}_l(h).
\end{equation}
Now we can state the following proposition, which immediately leads to Theorem \ref{PaperThmB}

\tsp
\newtheorem{pheene}[p1]{Proposition}
\begin{pheene}
Fix $\mbb\in(\R_{>0})^\n$ and set $V=\TmTr$. The restriction
\[
\Psibo=\Psi\restr{V\setminus V_0}\colon V\setminus V_0\to\mML\times\mML,
\]
where $V_0$ denotes the zero section of $V$, is bijective onto $\mFMLu\setminus\{(\O,\O)\}$.
\begin{proof}
Put
\[
\Psib=\Psi\restr{V}\colon V\to\mML\times\mML.
\]
Clearly $\Psib(V_0)=\{\vv\}$.\\
In the first place, let us show that $\Psib(V)\subseteq\mMLu$. 
From \eqref{saens}, if $\Psib\hta=\lpm$ then $\e{\pl}+\e{\ml}=0$, . Using \eqref{M} we get
\[
0=\ddt_{|0}\bigg(\ell_{E^{t\pl}_l(h)}(\de_i)+\ell_{E^{t\ml}_l(h)}(\de_i)\bigg)=-m(\de_i,\pl)-m(\de_i,\ml)
\]
for $i=1,\ldots,\n$.\\
In order to see that $\Psibo$ is injective onto $\mFMLu\setminus\{\vv\}$, we need to prove that for every $\lpm\in\mML\times\mML$ there is a 1:1 correspondence between $(\Psib)^{-1}\lpm$ and $e^{\pl}_r\cap \e{\ml}$. If $\lpm=\Psib\hta $ then $e_r^{\pl}$ and $e_l^{\ml}$ meet over $h$ by construction. Conversely, if there exists $h\in\mT$ such that $e_r^{\pl}(h)=\e{\ml}(h)$, then $\lpm=\Psib(h,e_r^{\pl}(h))$.\\
Now, we showed in Section \ref{Sez3} for $\lpm\neq\vv$ that $\e{\pl}+\e{\ml}$ has exactly one critical point if $\lpm$ is a filling couple, zero critical points otherwise (Remark \ref{Emily}). Therefore, $\Psibo$ is injective with image $\mFMLu\setminus\{\vv\}$.
\end{proof}
\end{pheene}


\bibliographystyle{amsplain}\nocite{*}\bibliography{Bib_ref}


\end{document}